\RequirePackage{lineno}
\documentclass[11pt]{article}

\usepackage{hyperref}
\usepackage[T1]{fontenc}
\usepackage[latin1]{inputenc}
\usepackage{amsmath,amsthm,amssymb}
\usepackage{mathabx}
\usepackage{xcolor}
\usepackage{authblk}
\usepackage[margin=1in]{geometry}

\usepackage{mathrsfs}

\usepackage{graphicx}
\usepackage{amsbsy}
\usepackage{yfonts}
\usepackage[all]{xy}
\usepackage{pgfplots}

\providecommand{\U}[1]{\protect\rule{.1in}{.1in}}

\newtheorem{prop}{Proposition}[section]
\newtheorem{cor}[prop]{Corollary}

\newtheorem{lem}[prop]{Lemma}

\newtheorem{theo}[prop]{Theorem}

\newcommand{\vertiii}[1]{{\left\vert\kern-0.25ex\left\vert\kern-0.25ex\left\vert #1
    \right\vert\kern-0.25ex\right\vert\kern-0.25ex\right\vert}}
    
    \newcommand{\goodchi}{\protect\raisebox{2pt}{$\chi$}}

\def\tr{\mbox{\rm Tr}}

\newcommand{\EE}{\mathbb{E}}

\newcommand{\LL}{\mathbb{L}}

\newcommand{\RR}{\mathbb{R}}

\newcommand{\Ca}{ {\cal C }}
\newcommand{\Da}{ {\cal D }}
\newcommand{\La}{ {\cal L }}
\newcommand{\Na}{ {\cal N }}
\newcommand{\Ka}{ {\cal K }}

\newcommand{\Ra}{ {\cal R }}

\newcommand{\Fa}{ {\cal F }}

\newcommand{\Oa}{ {\cal O }}
\newcommand{\Ia}{ {\cal I }}
\newcommand{\Xa}{ {\cal X }}
\newcommand{\Ma}{ {\cal M }}

\newcommand{\Ja}{ {\cal J }}

\newcommand{\point}{\mbox{\LARGE .}}

\newcommand{\cqfd}{\hfill\blbx \\}
\def\blbx{\hbox{\vrule height 5pt width 5pt depth 0pt}\medskip}

\def \RR{\mathbb{R}}

\def \EE{\mathbb{E}}

\def \LL{\mathbb{L}}

\numberwithin{equation}{section}

\makeatletter
\DeclareRobustCommand\frownotimes{\mathbin{\mathpalette\frown@otimes\relax}}
\newcommand{\frown@otimes}[2]{
  \vbox{
    \ialign{##\cr
      \hidewidth$\m@th#1{}_\frown$\kern-\scriptspace\hidewidth\cr
      \noalign{\nointerlineskip\kern-.1pt}
      $\m@th#1\otimes$\cr
    }
  }
}
\makeatother

\hypersetup{
    colorlinks=true,
    linkcolor=blue,
    filecolor=magenta,      
    urlcolor=cyan,
    pdftitle={Overleaf Example},
    pdfpagemode=FullScreen,
    }

\begin{document}

\begin{center}
{\Large \textbf{On Time Uniform  Wong-Zakai Approximation Theorems}}

\vspace{0.5cm}

PIERRE DEL MORAL$^{1}$, SHULAN HU$^{2}$, AJAY JASRA$^{3}$, HAMZA RUZAYQAT$^{3}$  \& XINYU WANG$^{4}$

{\footnotesize 
$^{1}$Institut de Mathematiques de Bordeaux, 33405 Bordeaux, FR\\
$^{2}$School of Statistics and Mathematics, Zhongnan University of Economics and Law, CN. \\
$^{3}$Applied Mathematics and Computational Science Program, \\ Computer, Electrical and Mathematical Sciences and Engineering Division, \\ King Abdullah University of Science and Technology, Thuwal, 23955-6900, KSA.\\
$^{4}$Wenlan School of Business, Zhongnan University of Economics and Law, CN.
} \\
\footnotesize{E-Mail: \verb|pierre.del-moral@inria.fr, hu_shulan@zuel.edu.cn|,
\verb|ajay.jasra@kaust.edu.sa, hamza.ruzayqat@kaust.edu.sa, wang_xin_yu@zuel.edu.cn|}

\begin{abstract}
We consider the long time behavior of Wong-Zakai approximations of stochastic differential equations. These piecewise smooth diffusion approximations are of great importance in many areas, such as those with ordinary differential equations associated to random smooth fluctuations; e.g.~robust filtering problems. 
In many examples, the mean error estimate bounds that have been derived in the literature can
grow exponentially with respect to the time horizon. 
We show in a simple example that indeed mean error estimates do explode exponentially in the time parameter, i.e.~in that case a Wong-Zakai approximation is only useful 
for extremely short time intervals. Under spectral conditions, we present some quantitative time-uniform convergence theorems, i.e.~time-uniform mean error bounds, yielding what seems to be the first results of this type for Wong-Zakai diffusion approximations.\\ 
\noindent\textbf{Keywords:} Wong-Zakai approximations, Stratonovich stochastic integral, It\^o's integral, stochastic differential equation, Brownian motion.\\
\noindent\textbf{Mathematics Subject Classification:} Primary 60H35; secondary 65C30.
\end{abstract}

\end{center}

\section{Introduction}
 Let $\overline{X}_{t}(x)$ be the flow of an ordinary differential equation in $\RR^{r}$, $r\in \mathbb{N}$, with 
    \begin{equation}
    \label{diff-def-over-X}
   d \overline{X}_{t}(x)=b\left(\overline{X}_{t}(x)\right)~dt+\sigma\left(\overline{X}_{t}(x)\right)~
    d\overline{B}_t\quad \mbox{\rm with initial condition}\quad \overline{X}_0(x)=x,
    \end{equation}
where $\overline{B}_t$ is a piecewise smooth approximation of  an $\overline{r}$-dimensional Brownian motion $B_t$, the drift function 
 $x\in \RR^{r}\mapsto b(x)\in \RR^{r}$ and the diffusion function 
 $x\in \RR^{r}\mapsto\sigma(x)\in \RR^{r\times \overline{r}}$ are some regular functions. Under appropriate conditions, the flow $\overline{X}_t(x)$ is a 
piecewise smooth approximation of the stochastic flow $X_t(x)$, with initial condition $X_0(x)=x$, defined by
the Stratonovitch stochastic differential equation
    \begin{equation}
    \label{diff-def}
    d X_{t}(x)=b\left(X_{t}(x)\right)~dt+\sigma\left(X_{t}(x)\right)~
    d^SB_t.
    \end{equation}
  The usage of such approximation schemes includes robust filtering e.g.~\cite{clark} and numerical approximation of stochastic differential equations. In the latter case, the approximation
  permits the application of numerical methods for ordinary differential equation over standard methods for stochastic differential equations such as Euler-Maruyama which 
  are typically designed to converge to an Ito stochastic differential equation; see e.g.~\cite{rum}.
 
 The convergence analysis of these smooth diffusion approximations have started 
  in the seminal article by Wong and Zakai in the mid-1960s~\cite{wong-zakai} in the context of one dimensional models, followed in the beginning of the 1970s by the article by Stroock and Varadhan~\cite{stroock-varadhan} for multi-dimensional models. Since then, there is a vast literature concerning the convergence analysis of these Wong-Zakai approximations, see for instance~\cite{friz,kon,kurtz,shmatkov,suss,tward}, as well as~\cite{acqui,brzez,gyon,gyon-2,hairer,tessi} in the context of stochastic partial differential equations and Hilbert-space valued processes. This is a non-exhaustive list of the literature; we refer the reader to the review article~\cite{tward} as well as  the more recent articles~\cite{ganguly,sahani} and the references therein for a more comprehensive coverage of the subject.
  
Most of these studies provide asymptotic almost sure convergence or weak convergence of the laws of Wong-Zakai approximations on compact time intervals toward a Stratonovich stochastic differential or partial differential equation. From the practical perspective, the few convergence estimates that we found in the literature, (i.e.~mean error bounds),
 grow exponentially fast with respect to the time horizon. This often means that Wong-Zakai approximation is only useful 
for extremely short time intervals, if indeed the bounds are sharp.
In addition, as underlined in~\cite{sahani}, much of the literature on Wong-Zakai approximations relies on the uniformly boundedness of the functions $(b,\sigma)$, thus does not apply to simple linear diffusions. The main technical difficulty in the stochastic analysis of Wong-Zakai approximations is that it requires refined estimates of anticipative stochastic integrals. Combining truncation techniques with the rough-path theory pioneered by T. Lyons~\cite{terry}, new almost sure and asymptotic  $\LL_2$ convergence theorems in the unbounded case are provided in the recent article~\cite{sahani}.

\subsection{Contribution and Structure}

The main contribution of this article is to present non asymptotic and time uniform estimates for Wong-Zakai approximations. To the best of our knowledge these bounds are the first to appear in the literature. The main importance from the practical perspective is to provide guarantees on when such approximations can be accurate over long-time periods.

To understand the kind of quantitative estimates we might expect, we consider 
the simplest context of one dimensional models with a prescribed
 time mesh  $t_n< t_{n+1}$ with $\epsilon:=(t_{n+1}-t_n)$. We let
 $  \overline{B}_t$ be the piecewise linear smooth anticipative one dimensional process  defined on the time intervals $t_n\leq t< t_{n+1}$ by
\begin{eqnarray}
\overline{B}_t&:=&B_{t_n}+\frac{t-t_n}{\epsilon}~(B_{t_{n+1}}-B_{t_n})\label{ADT}\\
\Longrightarrow B_t-\overline{B}_t&=&\frac{t_{n+1}-t}{\epsilon}~(B_t-B_{t_n})+\frac{t-t_n}{\epsilon}~(B_t-B_{t_{n+1}})=\sqrt{\epsilon~\frac{(t_{n+1}-t)}{\epsilon}~\frac{(t-t_n)}{\epsilon}}~\goodchi_t \nonumber
\end{eqnarray}
where $\goodchi_t$ is a centered Gaussian random variable with unit variance.
This yields, for any $p\geq 1$, the even moments formula
 $$
\EE\left((B_t-
\overline{B}_t)^{2p}\right)=\frac{(2p)!}{p!2^p} ~\epsilon^p~\left(\frac{(t_{n+1}-t)}{\epsilon}~\frac{(t-t_n)}{\epsilon}\right)^p,
$$
from which we find the time uniform estimates 
 \begin{equation}\label{Hyp-over-B}
\sup_{t\geq 0}\EE\left(\vert B_t-
\overline{B}_t\vert^{q}\right)^{1/q} 
\leq c_q~\sqrt{\epsilon},
\end{equation}
for all $q\geq 1$.
The above estimate is also satisfied for  the smooth approximation
 \begin{equation}\label{IOU}
\overline{B}_t:=\int_0^t~Y_s~ds
\quad \mbox{\rm with}\quad 
dY_t=-\frac{1}{\epsilon}~Y_t~dt+\frac{1}{\epsilon}~dB_t\quad \mbox{\rm and}\quad 
Y_0=0.
\end{equation}
Note that
$$
B_t-
\overline{B}_t=\epsilon Y_t=\int_0^t~e^{-(t-u)/\epsilon}~dB_u=  \sqrt{\frac{\epsilon}{2}~(1-e^{-2t/\epsilon})}~\goodchi_t.
$$
This yields, for any $p\geq 1$, the formula
 $$
\EE\left((B_t-
\overline{B}_t)^{2p}\right)=\frac{(2p)!}{p!2^{2p}}~\epsilon^p~(1-e^{-2t/\epsilon})^p\leq \frac{(2p)!}{p!2^{2p}}~\epsilon^p.
 $$

Consider the one dimensional  version of the stochastic flows $X_t(x)$ and $\overline{X}_t(x)$ defined in \eqref{diff-def-over-X} and \eqref{diff-def} with the piecewise linear approximation discussed in \eqref{ADT} and 
 \begin{equation}\label{LG-ref}
b(x)=a~x\quad\mbox{\rm for some}\quad a\in \RR\quad \mbox{\rm with unit diffusion function}\quad \sigma(x)=1.
\end{equation}
In this context, two scenarios arise:
\begin{itemize}
\item When $a<0$, we have the time uniform estimate
 \begin{equation}\label{LG-stable}
 \sup_{n\geq 0}\EE\left((X_{t_n}(x)-\overline{X}_{t_n}(x))^2\right)\leq c(\epsilon)~\epsilon^2
\quad \mbox{with}\quad c(\epsilon):=~e^{2\vert a\vert \epsilon}~{\vert a\vert}/{24}.
 \end{equation}
\item When $a>0$, we have the exponential growth estimates
  \begin{equation}\label{LG-unstable}
\epsilon^2~c_-(\epsilon)~\left(e^{2at_n}-1\right)\leq 
\EE\left((X_{t_n}(x)-\overline{X}_{t_n}(x))^2\right)\leq \epsilon^2~c_+(\epsilon)~\left(e^{2at_n}-1\right)
   \end{equation}
with the parameters
$$
c_-(\epsilon):=\frac{a}{24}~e^{-2a\epsilon}\quad \mbox{and}\quad
c_+(\epsilon):=\frac{a}{24}~e^{2a\epsilon}.
$$
\end{itemize}
The detailed proof  of \eqref{LG-stable} and \eqref{LG-unstable} is provided in the Appendix. Choosing $a=2$ and $t_n\geq 10$ in the bounds in \eqref{LG-unstable} yields that
$e^{2at_n}=e^{40}$ and  is  $16$ times larger than the estimated age of the universe given by $13.7~e^{28}$ years. These exponential bounds are not only impractical from a numerical perspective, they also suggest that Wong-Zakai approximations can only be used for stable-like stochastic  processes.
 In this context, classical perturbation bounds combining Lipschitz type inequalities with Gronwall lemma also yield exceedingly pessimistic global estimates that grows exponentially fast with respect to (w.r.t.) the time horizon.

We also seek to obtain useful quantitative and time uniform estimates which are valid under a single set of easily checked conditions that only depend on the parameters of the model. Throughout the article we assume that  the diffusion  function  $\sigma$ as well as its first and second derivatives  are uniformly bounded. In addition,
$ x\mapsto b(x)$ have continuous and uniformly bounded derivatives up to the second order.  Various techniques presented in the article and many results can be separately and readily extended to more general models with weaker and abstract assumptions.

This article is structured as follows. In Section \ref{sec:not_res} we provide the notation used in this article and state some of the main results. In Section~\ref{sec:prel_res}, we present some preliminary results. Sections~\ref{theo-sigma-ct-proof-i}, \ref{sec:proof_theo_loc} and \ref{theo-unif-int-strato-proof-1} are devoted for the proofs of the main theorems which are to be stated.
In particular, Section \ref{theo-sigma-ct-proof-i} contains the proof of Theorem \ref{theo-sigma-ct}, Section \ref{sec:proof_theo_loc} has the proofs of Theorems \ref{theo-loc-intro} and \ref{theo-unif-2} and Section \ref{theo-unif-int-strato-proof-1}  houses the proof of Theorem \ref{theo-unif-int-strato}. Some technical results are proved in the Appendix.

\section{Notation and Main Results}\label{sec:not_res}

\subsection{Basic Notation}
With a slight abuse of notation, we denote by $I$ the identity $(r\times r)$-matrix, for any $r \in \mathbb{N}$. We also denote by $\Vert\point\Vert$ any (equivalent) norm on a finite dimensional vector space over $\RR$. All vectors are column vectors by default. Next, we introduce some matrix notation needed from the onset. We denote by $\tr(A)$, $\Vert A\Vert_{2}:=\lambda_{\max}(AA^{\prime})^{1/2}=\lambda_{\max}(A^{\prime}A)^{1/2}$, $\Vert A\Vert_{\tiny Frob}:=\tr(AA^{\prime})^{1/2}$ and $\rho(A):=\lambda_{\max}((A+A^{\prime})/2)$, respectively, the trace, the spectral norm, the Frobenius norm, and the logarithmic norm of some matrix $A$. We denote by $A^{\prime}$ the the transpose of $A$ and $\lambda_{\max}(\point)$ the largest eigenvalue. The spectral norm is sub-multiplicative or $\Vert A B\Vert_{2}\leq \Vert A\Vert_{2} \Vert B\Vert_{2}$ and compatible with the Euclidean norm for vectors; by that we mean that for a vector $x$, we have $\Vert A x \Vert \leq \Vert A\Vert_{2} \Vert x\Vert$. Let $[n]$ be the set of $n$ multiple indexes $i=(i_1,\ldots,i_n)\in \Ia^n$ over some finite set $\Ia$. We denote by $(A_{i,j})_{(i,j)\in [p]\times [q]}$ the  entries of a $(p,q)$-tensor $A$ with index set $\Ia$ for $[p]$ and $\mathcal{J}$ for $[q]$. For the sake of brevity, the index sets will be implicitly defined through the context. For a given $(p_1,q)$-tensor $A$ and a given $(q,p_2)$ tensor $B$, $AB$ and $B^{\prime}$ are a $(p_1,p_2)$-tensor and a $(p_2,q)$-tensor, respectively, with entries
given by
\begin{equation}\label{tensor-notation}
\forall (i,j)\in [p_1]\times [p_2],\qquad
(AB)_{i,j}=\sum_{k\in [q]}A_{i,k}~B_{k,j}\quad \mbox{\rm and}\quad B_{j,k}^{\prime}:=B_{k,j}.
\end{equation}
The symmetric part of a $(p,p)$-tensor is the $(p,p)$-tensor $A_{\tiny sym}$ with entries
 $$
 \forall (i,j)\in [p]\times [p],\qquad(A_{\tiny sym})_{i,j}=(A_{i,j}+A_{j,i})/2.
 $$
We consider the Frobenius inner product given for any $(p,q)$-tensors $A$ and $B$ by
$$
\langle A,B\rangle_F:=\tr(AB^{\prime})=\sum_{i}(AB^{\prime})_{i,i}.
$$
For any $(p,q)$-tensors $A$ and $B$, we also check the Cauchy-Schwartz inequality
$$
\langle A,B\rangle_{\tiny Frob}^2\leq \Vert A\Vert_{\tiny Frob}~\Vert B\Vert_{\tiny Frob}.
$$
For any tensors $A_1,A_2$ with appropriate dimensions we have the inequality
$$
\Vert A_1A_2\Vert_{\tiny Frob}\leq \Vert A_1\Vert_{\tiny Frob}~\Vert A_2\Vert_{\tiny Frob}.
$$
Given a column vector-valued smooth function $f=(f^i)_{1\leq i\leq p}$
from $\RR^{r}$ into $\RR^{p}$
we denote by $\nabla f$ the gradient $(r,p)$-matrix defined by 
\begin{equation}
\label{grad-def}
\nabla f:=\left[\nabla f^1,\ldots,\nabla f^p\right]\quad \mbox{\rm with the gradient column vector}\quad\nabla f^i=\left[\begin{array}{c}
\partial_{1}f^i\\
\vdots\\
\partial_{r}f^i
\end{array}\right].
\end{equation} 
In the above display, $\partial_{j}f^i$, with $1\leq j\leq r$ and $1\leq i\leq p$, stands for the derivative of the function component
$x=(x_1,\ldots,x_r)\mapsto f^i(x)$ with respect to the $j$-th coordinate. In what follows, we also use the notation $\partial_{j,k}=\partial_j\partial_k=\partial_k\partial_j$, with $1\leq j,k\leq r$.
We also consider the following for $(2,2)$-tensors
\begin{eqnarray}
\left[\nabla f(x)\otimes \nabla f(x)\right]_{(i,j),(k,l)}&:=&\nabla f(x)_{i,k}\nabla f(x)_{j,l}=\left[\nabla f(x)\otimes \nabla f(x)\right]^{\prime}_{(k,l),(i,j)},\nonumber
\end{eqnarray}
where $\otimes$ denotes the tensor product. The notation $\partial f$ is {\em only} used in the context of one-dimensional models, $\partial^n f$ is the $n$-th derivative of a function $f(x)$ from $\RR$ into itself. We remark that the usage of vector-type gradient and matrix-type Hessian notation for one-dimensional models is potentially confusing. For any smooth column vector-valued function $g=(g^i)_{1\leq i\leq r}$ from $\RR^{r}$ into itself, we have
\begin{align*}
&h:=f\circ g\\
\\
\Longrightarrow~
&h^i(x):=(f^i\circ g)(x):=f^i(g(x))\\
\\
\Longrightarrow~ &\nabla h(x)_{k,i}=\partial_{k}h^i(x)=\sum_{1\leq l\leq r}(\partial_{l}f^i)(g(x))~\partial_{k}g^l(x)=\sum_{1\leq l\leq r} ~\nabla g(x)_{k,l}~(\nabla f)(g(x))_{l,i}.
\end{align*}
Equivalently, we have the composition rule matrix formula
$$
\nabla (f\circ g)(x)=\nabla g(x)~(\nabla f)(g(x)).
$$
Building on this notation, let
\begin{equation}\label{Hessian-def}
\nabla^2 f=\left[\nabla^2 f^1,\ldots,\nabla^2 f^p\right]\quad \mbox{\rm with the $(r\times r)$-matrix}\quad\nabla^2 f^i=\left[\begin{array}{ccc}
\partial_{1,1}f^i&\ldots&\partial_{1,r}f^i\\
\vdots&\ldots&\vdots\\
\partial_{r,1}f^i&\ldots&\partial_{r,r}f^i\
\end{array}\right].
\end{equation} 
Thus, the Hessian $\nabla^2 f$ and its transpose $(\nabla^2 f)^{\prime}$  are defined by the $(2,1)$-tensor and $(1,2)$-tensor with entries
  $$
(\nabla^2 f)_{(i,j),k}=(\nabla^2 f^k)_{i,j}=\partial_{i,j}f^k= (\nabla^2 f)^{\prime}_{k,(i,j)}.
 $$
In the same vein, we have
\begin{align*}
&(\nabla^2 (f\circ g))(x)_{(i,j),k}\\
&=\partial_{i,j}\left(f^k(g(x))\right)\\
&=\sum_{1\leq l,m\leq r}(\partial_{l,m}f^k)(g(x))~\partial_{i}g^l(x)\partial_{j}g^m(x)+
\sum_{1\leq l\leq r}(\partial_{l}f^k)(g(x))~\partial_{i,j}g^l(x)\\
&=\sum_{1\leq l,m\leq r} 
\left[ \nabla g(x)\otimes \nabla g(x)\right]_{(i,j),(l,m)}~(\nabla^2 f)_{(l,m),k}+
\sum_{1\leq l,m\leq r} (\nabla^2 g)_{(i,j),l}~(\nabla f)(g(x))_{l,k}.
\end{align*}
Equivalently, we have the second order composition rule tensor formula
$$
\nabla^2 (f\circ g)(x)=\left[ \nabla g(x)\otimes \nabla g(x)\right](\nabla^2f)(g(x))+\nabla^2 g(x)~(\nabla f)(g(x)).
$$
This notation allows one to compactly represent the second order term of the Taylor expansion of the the vector valued function $f$. For a column vector $y=(y_1,\ldots,y_r)'$, we have
 \[
 \left [
 \begin{array}{c}
 y^{\prime}~\nabla^2 f^1(x)~y \\
 \vdots \\
 y^{\prime}~\nabla^2 f^p(x)~y \\ 
 \end{array}
 \right ] = (\nabla^2 f(x))^{\prime}~yy^{\prime}\Longleftrightarrow
\left(  (\nabla^2 f(x))^{\prime}~yy^{\prime}\right)_{k,1}=\sum_{ (i,j)}
(\nabla^2 f(x))^{\prime}_{k,(i,j)}~(yy')_{(i,j),1},
 \]
 where we have regarded the matrix $a:=yy'$ as the $(2,1)$-tensor $a^{(k,l)}=(yy')_{(k,l),1}=y_ky_l$, and the summation in the above display is taken over all $2$-indices $(i,j)\in \{1,\ldots,r\}^2$.
For any 
 $r$-column vector valued function $b(x)$ with entries $b^i(x)$ and for any  $(r\times r)$-matrix function $x\in\RR^r\mapsto a(x)=(a(x)^{i,j})_{1\leq i,j\leq r}$, we have
  $$
 \begin{array}{rcl}
 \displaystyle\left( \nabla f(x)^{\prime}~b(x)\right)^k&=& \displaystyle\sum_{1\leq i\leq r}( \nabla f(x))^{\prime}_{k,i}~b^i(x)=\sum_{1\leq i\leq r}~\partial_{i} h^k(x)~b^i(x):=\langle \nabla f^k(x),b(x)\rangle\\
 \\
 \displaystyle\left(\nabla^2 f(x)^{\prime}~a(x)\right)^k&=& \displaystyle\sum_{1\leq i,j\leq r}( \nabla^2 f(x))^{\prime}_{k,(i,j)}~a(x)^{i,j}\\
 &&\\
 &=& \displaystyle\sum_{1\leq i,j\leq r}~ \partial_{i,j} f^k(x)~a(x)^{i,j}:=\langle \nabla^2f^k(x),a(x)\rangle_{\tiny Frob}~.
\end{array}
 $$
 In a more compact form, the above formula takes the form
  \begin{equation}\label{tensor-product-2-2}
  \nabla f(x)^{\prime}~b(x)=\left[\begin{array}{c}
\langle\nabla f^1(x),b(x)\rangle\\
\vdots\\
\langle\nabla f^p(x),b(x)\rangle
\end{array}\right]\quad\mbox{\rm and}\quad
  \nabla^2 f(x)^{\prime}~a(x)=\left[\begin{array}{c}
\langle\nabla^2f^1(x),a(x)\rangle_{\tiny Frob}\\
\vdots\\
\langle\nabla^2f^p(x),a(x)\rangle_{\tiny Frob}
\end{array}\right].
\end{equation} 

  Throughout, unless otherwise is stated we write $c$ for some positive constants whose values may vary from line to line, and we write $c_{\alpha}$, as well as $c(\beta)$ and $c_{\alpha}(\beta)$ when their values may depend on  some parameters $\alpha,\beta$ defined on some parameter sets. We also set $a\wedge b=\min(a,b)$, $a\vee b=\max(a,b)$, and $a_+=a\vee 0$ for $a,b\in \RR$.

\subsection{Statement of Some Main Results}
Let $r,\overline{r}\geq 1$ be some integer  parameters. Consider an $\overline{r}$-dimensional standard Brownian process $B_t$ starting at the origin, that is, $B_0=0$.  As a rule, $B_t$ is represented by an $\overline{r}$-column vector with $i$-entries $B^i_t$, for $1\leq i\leq \overline{r}$. The polygonal approximation $\overline{B}_t$ of $B_t$ on some time mesh  $t_n< t_{n+1}$ is defined in a vector form by the formula \eqref{ADT}. In what follows,
$\overline{B}_t$ stands for some piecewise smooth approximation of  an $r$-dimensional Brownian motion $B_t$ such that, for any $p\geq 1$, we have the time uniform estimates
 \begin{equation}\label{Hyp-over-B}
\sup_{t\geq 0}\EE\left(\Vert B_t-
\overline{B}_t\Vert^{p}\right)^{1/p} 
\leq c_p~\sqrt{\epsilon}.
\end{equation}
State variables  of the differential equations \eqref{diff-def-over-X} and \eqref{diff-def} are represented by $r$-column vectors with $i$-entries $\overline{X}_t^{\,i}(x)$ and  $X_t^i(x)$, for $1\leq i\leq r$.

Consider the following spectral condition:
\begin{quote}
{\em
\noindent $(H_b)$ The log-norm $\rho\left(\nabla b(x)\right)$ of the Jacobian matrix $\nabla b(x)$ of the drift function $b$ is uniformly negative, that is we have
\begin{equation}\label{hyp-refa} 
-\lambda_b:=\sup_{x\in \RR^r}\rho\left(\nabla b(x)\right)<0.
\end{equation}
}
\end{quote}
This condition is clearly met for the stable linear Gaussian models discussed in \eqref{LG-stable}. 
As shown in Theorem~\ref{theo-moments}, the above condition allows one to control, uniformly in time, the absolute moments of the smooth diffusion approximation associated with  the  polygonal approximation  $  \overline{B}_t$ defined in \eqref{ADT}; that is, we have that
 \begin{equation}\label{unif-intro}
 (H_b)\Longrightarrow \forall p\geq 1,\quad \sup_{t\geq 0}\EE\left(\Vert \overline{X}_t(x)\Vert^p\right)\leq c_p~(1+\Vert x\Vert).
 \end{equation}
 
 We are now in a position to state our first result. Below we will say that $\sigma$ is a homogeneous diffusion matrix, by which we mean that $\sigma(x)=\sigma(0)$.
 The proof of the following theorem can be found in Section \ref{theo-sigma-ct-proof-i}.

\begin{theo}\label{theo-sigma-ct}
Assume $(H_b)$, $r=\overline{r}$ and that $\sigma$ is a homogeneous diffusion matrix. Then there exists $\lambda>0$, such that for any $t\geq 0$, we have the almost sure uniform estimate
$$
\Vert X_{t}-\overline{X}_{t}\Vert_{\infty}:=
\sup_{x\in\RR^r}\Vert X_{t}(x)-\overline{X}_{t}(x)\Vert \leq c~\left(\Vert B_t-
\overline{B}_t\Vert+\left(\int_0^t~e^{-2\lambda(t-s)}~
\Vert B_s-
\overline{B}_s\Vert^2~
ds\right)^{1/2}\right).
$$
In addition, the mean error estimates \eqref{Hyp-over-B} yields, for any $p\geq 1$, the space-time uniform estimate
\begin{equation}\label{unif-intro-X-over}
\sup_{t\geq 0}\EE\left(\Vert X_{t}-\overline{X}_{t}\Vert_{\infty}^{p}\right)^{1/p}\leq c_p~\sqrt{\epsilon}.
\end{equation}
\end{theo}

Providing Wong-Zakai uniform estimates for diffusion flows equipped with non homogeneous diffusion matrices is technically more involved. As expected, these stochastic models can be handled when a transformation to unit diffusion exists. For instance, consider a one dimensional diffusion
with $0<\sigma_-\leq \sigma(x)\leq \sigma_+$ and the Lamperti transformation 
$$
\Xa_t:=\theta(X_t)\quad \mbox{\rm and}\quad
\overline{\Xa}_t:=\theta(\overline{X}_t)\quad \mbox{\rm with}\quad
\theta(x):=\int_{0}^x~\frac{1}{\sigma(y)}~dy.
$$
In this context, we have
\begin{equation}\label{lamperti}
d\Xa_t=b^{\theta}(\Xa_t)dt+dB_t\quad \mbox{\rm and}\quad
d\overline{\Xa}_t=b^{\theta}(\overline{\Xa}_t)dt+d\overline{B}_t,
\end{equation}
with the drift function
$$
b^{\theta}(x):=b(\theta^{-1}(x))/\sigma(\theta^{-1}(x)).
$$
The proof of \eqref{lamperti} follows well known arguments, however, is provided in the Appendix on page~\pageref{lamperti-proof}.
In this context, using the fact that
$$
\vert x_2-x_1\vert/\sigma_+\leq \vert \theta(x_1)-\theta(x_2)\vert \leq  \vert x_2-x_1\vert/\sigma_-~,
$$
without further work, we readily check the following corollary.
\begin{cor}
 Assume that the log-norm $\rho\left(\nabla b^{\theta}(x)\right)$ associated with the one dimensional model (\ref{lamperti}) is uniformly negative. In this situation, the mean error estimates (\ref{Hyp-over-B}) yields, for any $p\geq 1$, the space-time uniform estimate
$$
\sup_{t\geq 0}\EE\left(\Vert X_{t}-\overline{X}_{t}\Vert_{\infty}^{p}\right)^{1/p}\leq c_p~\sqrt{\epsilon}.
$$
\end{cor}

To handle more general situation, recall that
the Stratonovitch stochastic differential equation (\ref{diff-def}) can be rewritten in It\^o's form 
as follows
    \begin{equation}\label{diff-def-Ito}
    d X_{t}(x)=b_{\sigma}\left(X_{t}(x)\right)~dt+\sigma\left(X_{t}(x)\right)~
    dB_t ,   
    \end{equation} 
with the drift function
$$
b_{\sigma}^i:=b^i+\frac{1}{2}~\sum_{1\leq j\leq \overline{r}}
\sum_{1\leq k\leq r}~\sigma^k_{j}~\partial_k \sigma^i_{j}~.
$$
Consider the following spectral condition:
\begin{quote}
\noindent{\em
 $(H_{b_{\sigma}})$ 
 The log-norm $\rho\left(\nabla b_{\sigma}(x)\right)$ of the Jacobian matrix $\nabla b_{\sigma}(x)$ of the drift function $b_{\sigma}$ is uniformly negative, that is we have
\begin{equation}\label{hyp-refa-sig} 
-\lambda_{b_{\sigma}}:=\sup_{x\in \RR^r}\rho\left(\nabla b_{\sigma}(x)\right)<0.
\end{equation}
}
\end{quote}
This condition ensures that, for any $p\geq 1$ and $x\in\RR^r$, we have
the uniform absolute moments estimates
\begin{equation}\label{hyp-mom}
(H_{b_{\sigma}})\Longrightarrow
\sup_{t\geq 0}\EE\left(\Vert X_{t}(x)\Vert^p \right)^{1/p}\leq c_p~(1+\Vert x\Vert).
\end{equation}
A proof of \eqref{hyp-mom} is provided in the Appendix on page~\pageref{hyp-mom-proof} (see also  Lemma 2.1 in~\cite{dmsumeet}). Note that
$$
\eqref{hyp-mom}\Longrightarrow
 \EE\left(\Vert X_{t}(x)-x\Vert^p \right)^{1/p}\leq c_p~\left(
 t~(1+\Vert x\Vert)+\sqrt{t}\right).
$$

We further assume that $\overline{B}_t$ is the polygonal approximation  of $B_t$ on some time mesh  $t_n< t_{n+1}$ defined in a vector form by the formula (\ref{ADT}). 
We denote by $\Fa_t$ the $\sigma$-field generated by the Brownian motion $B_s$ with $s<t$, and for any $n\geq 1$, we let
$
\goodchi_{t_n}
$ be the $r$-column vector with $i$-th entries $$\goodchi^i_{t_n}:=(B^i_{t_n}-B^i_{t_{n-1}})/\sqrt{\epsilon}.
$$
We also set
$$
\Delta X_{t_{n}}:= X_{t_{n}}-X_{t_{n-1}}\quad \mbox{\rm and}\quad
\Delta \overline{X}_{t_{n}}:= \overline{X}_{t_{n}}- \overline{X}_{t_{n-1}}.
$$
In this notation, we have the following Milstein-type decompositions which is proved in Section~\ref{sec:proof_theo_loc}.
\begin{theo}\label{theo-loc-intro}
Assume that \eqref{unif-intro}
and \eqref{hyp-mom} hold and the following commutation property is satisfied
    \begin{equation}\label{ref-commutation}
\nabla \sigma_k(x)^{\prime} \sigma_{l}(x)=
\nabla \sigma_l(x)^{\prime} \sigma_{k}(x)
 \end{equation} 
for any $1\leq k,l\leq \overline{r}$ and $x\in\RR^r$.
Then for any $n\geq 1$, we have the second order Milstein-type decompositions
\begin{eqnarray}
\Delta  \overline{X}_{t_{n}}&=&b_{\sigma}\left( \overline{X}_{t_{n-1}}\right)~\epsilon+\sigma\left( \overline{X}_{t_{n-1}}\right)\sqrt{\epsilon}~\goodchi_{t_{n}}\label{delta-overX-intro} \\
&&\hskip2.2cm\displaystyle+~\frac{1}{2}~\sum_{1\leq k,l\leq \overline{r}}
\nabla_{\sigma_l} \sigma_k(\overline{X}_{t_{n-1}})~\left(\goodchi_{t_{n}}^k~\goodchi_{t_{n}}^l-1_{k=l}\right)~\epsilon +\epsilon^{3/2}~\Delta \overline{R}_{t_{n}}~,\nonumber
\end{eqnarray}
as well as the second order expansion
\begin{eqnarray}
 \Delta X_{t_{n}}&=&b_{\sigma}\left(X_{t_{n-1}}\right)~\epsilon+\sigma\left(X_{t_{n-1}}\right)\sqrt{\epsilon}~\chi_{t_{n}} +\epsilon~\Delta M_{t_{n}}\label{delta-X-intro} \\
&&\hskip2.2cm\displaystyle+~\frac{1}{2}~\sum_{1\leq k,l\leq \overline{r}}
\nabla_{\sigma_l}\sigma_k(X_{t_{n-1}})\left(\goodchi_{t_{n}}^l\goodchi_{t_{n}}^k-1_{k=l}\right)~\epsilon+\epsilon^{3/2}~\Delta R_{t_{n}}~.\nonumber
\end{eqnarray}
In the above display  $\Delta M_{t_{n}}$ stands for a martingale increment and $\left(\Delta \overline{R}_{t_{n}},\Delta R_{t_{n}}\right)$ are remainder terms such that, for any  $x\in\RR^r$, $n\geq 1$ and $p\geq 1$, we have
the time uniform almost sure estimates
\begin{equation}\label{est-r-X}
\EE\left(\Vert \Delta M_{t_{n}} (x)\Vert^{p}~|~\Fa_{t_{n-1}}\right)^{1/p}\vee\EE\left(\Vert \Delta R_{t_{n}} (x)\Vert^{p}\right)^{1/p}\vee\EE\left(\Vert \Delta \overline{R}_{t_{n}} (x)\Vert^{p}\right)^{1/p}\leq c_{p}~(1+\Vert x\Vert).
\end{equation}
In addition, for any $1\leq i,j,k\leq r$,  we have
\begin{equation}\label{est-r-M}
\EE\left( \goodchi^i_{t_{n}}\Delta M^j_{t_{n}} (x)~|~\Fa_{t_{n-1}}\right)=0=\EE\left( \goodchi^i_{t_{n}}\goodchi^j_{t_{n}}~\Delta M^k_{t_{n}} (x)~|~\Fa_{t_{n-1}}\right).
\end{equation}
 \end{theo}
The next condition is a strengthening of the spectral condition  $(H_{b_{\sigma}})$:
\begin{quote}
\noindent{\em
 $(H_{\sigma})$ There exists some $\lambda_{\sigma}$ such that, for any $x\in \RR^r$, we have
$$
S_{\sigma}(x):=\left(\nabla b_{\sigma}(x)\right)_{\tiny sym}+2^{-1}\sum_{1\leq k \leq r}\nabla \sigma_{k}(x)(\nabla\sigma_{k}(x))^{\prime}\leq
-\lambda_{\sigma}~I .
$$}
\end{quote}
Condition  $(H_{\sigma})$ is reminiscent of the $\LL_2$-norm control of the tangent process $\nabla X_{t}(x)$ w.r.t the Frobenius norm. More precisely, as shown in Section 3 in~\cite{dmsumeet}, we have
$$
 \begin{array}{l}
\displaystyle d  \,\Vert \nabla X_{t}(x)\Vert^2_{\tiny Frob}
=2~\tr\left[\nabla X_{t}(x) ~S_{\sigma}\left(X_{t}(x)\right)~\nabla X_{t}(x)^{\prime}\right]~dt+d \Ma_{t}(x),
\end{array} 
$$
with the martingale
$$
d \Ma_{t}(x)=\sum_{1\leq k \leq r} \tr\left(\nabla X_{t}(x) \left[\nabla\sigma_{k}\left(X_{t}(x)\right)+\nabla\sigma_{k}\left(X_{t}(x)\right)^{\prime}\right]\nabla X_{t}(x)^{\prime}\right)~dB^k_t.
$$
This yields the time uniform $\LL_2$-norm estimates
$$
(H_{\sigma})\Longrightarrow \EE\left(\Vert \nabla X_{t}(x)\Vert^2_{\tiny Frob}\right)^{1/2}~\leq~\sqrt{r}~ e^{-\lambda_{\sigma}t}.
$$
The estimate of higher absolute moments requires more stringent conditions.  We have the following result whose proof is in Section~\ref{theo-unif-2-proof}. 
\begin{theo}\label{theo-unif-2}
Assume  $(H_{\sigma})$ and that \eqref{ref-commutation} and \eqref{unif-intro} hold. Then we have the uniform variance estimate
\begin{equation}\label{unif-X-n}
\sup_{n\geq 0}\EE\left(\Vert X_{t_{n}}(x)-\overline{X}_{t_{n}}(x)\Vert^2\right)^{1/2}\leq c~
\epsilon^{1/4}~(1+\Vert x\Vert).
\end{equation}
\end{theo}

We end this section with a refined and extended multidimensional version of \cite[Theorem 6.1]{brzez-0}. Let $(s,x)\in [0,t]\times \RR^r\mapsto F_{s,t}(x)\in \RR^{\overline{r}}$ be a twice differentiable column vector-valued  function. Assume there exists some $\lambda\geq 0$ such that, for any $0\leq u\leq s\leq t$ and $x\in \RR^r$, we have
\begin{equation}\label{hyp-Fst}
\left\{\begin{array}{l}
\Vert  \nabla F_{s,t}(x)\Vert\vee  \Vert \partial_s\nabla F_{s,t}(x)\Vert \vee\Vert  \nabla^2 F_{s,t}(x)\Vert\leq c~e^{-\lambda(t-s)}\\
\\
\mbox{\rm as well as}\quad \Vert \partial_sF_{s,t}(x)\Vert
 \leq c~e^{-\lambda(t-s)}~(1+\Vert x\Vert)\quad\mbox{\rm and}\quad \Vert F_{s,t}(x)\Vert \leq c~(1+\Vert x\Vert)~.
\end{array}\right.
\end{equation}
For any $\lambda>0$ we  set
$$
e_{\lambda}(t):=\int_0^{t}~e^{-\lambda s}~ds=\frac{1}{\lambda}\left(1-e^{-\lambda t}\right)\quad\mbox{\rm and}\quad e_{0}(t):= t \quad\mbox{\rm for $\lambda=0$.}
$$
We have the following result whose proof is proved in Section~\ref{theo-unif-int-strato-proof-1}.
\begin{theo}\label{theo-unif-int-strato}
Assume $(H_{b_{\sigma}})$ and \eqref{hyp-Fst} holds for some $\lambda\geq 0$. Then
for any time horizon $t\geq 0$, we have the uniform Wong-Zakai estimate
\begin{equation}\label{over-B-strato}
\begin{array}{l}
\displaystyle\int_{0}^{t}~F_{s,t}(X_{s}(x))^{\prime}~d\overline{B}_s
=\int_{0}^{t}~F_{s,t}(X_{s}(x))^{\prime}~dB_s+\frac{1}{2}
\int_{0}^{t}~\tr\left(\nabla_{\sigma} F_{s,t}(X_{s}(x))\right)~ds+R_{t}(F)(x)~,
\end{array}
\end{equation}
with the remainder term $R_{t}(F)(x)$
such that
$$
\EE\left( R_t(F)(x)^2\right)^{1/2}\leq 
c~\sqrt{\epsilon}~\left(1+\left(e_{\lambda}(t)\vee e_{2\lambda}(t)^{1/2}\right)\right)~(1+\Vert x\Vert).
$$
When \eqref{hyp-Fst} is met for some $\lambda>0$, we also have the uniform estimate
$$
\sup_{t\geq 0}\EE\left(R_t(F)(x)^2\right)^{1/2}\leq c~\sqrt{\epsilon}~(1+\Vert x\Vert).
$$
\end{theo}

Consider one dimensional diffusions
with $0<\sigma_-\leq \sigma(x)\leq \sigma_+$. Following the proof of \cite[Theorem 1a]{wong-zakai}, the uniform estimates stated in \eqref{unif-intro-X-over} ensure that
\begin{equation}\label{fin-s}
\int_0^t f(\overline{X}_s(x))~d\overline{B}_s
=\int_0^t f(X_s(x))~dB_s+\frac{1}{2}\int_0^t \sigma(X_s(x))~\partial f(X_s(x))~ds+\overline{R}_t(f)(x),
\end{equation}
for some remainder function $\overline{R}_t(f)(x)$, such that
$$
\EE\left(\overline{R}_t(f)(x)^2\right)^{1/2}\leq c~(1+t)~\sqrt{\epsilon}.
$$
The proof follows word-for-word the proof of \cite[Theorem 1a]{wong-zakai}, however, it is still provided in the Appendix on page~\pageref{fin-s-proof}.


\section{Some Preliminary Results}
\label{sec:prel_res}  

 The drift function $b(x)$ is a
 $r$-column vector valued function with entries $b^i(x)$ and $\sigma(x)$ is a $(r\times \overline{r})$-matrix valued function with $(i,j)$-entries $\sigma_{j}^i(x)$, for any $1\leq i\leq r$
 and $1\leq j\leq \overline{r}$. In this notation, we have
 $$
 \sigma=\left[ \sigma_1,\ldots, \sigma_{\overline{r}}\right]
 \quad \mbox{\rm with column vectors}\quad  \sigma_j=\left[\begin{array}{c}
 \sigma_j^1\\ 
\vdots\\
 \sigma_j^r
\end{array}\right].
 $$
For any column vector-valued
 smooth function $f=(f^i)_{1\leq i\leq p}$ we denote by $ \nabla_{\sigma}f(x)
 $ the $(\overline{r},p)$-matrix valued function with entries $(j,i)$-th entries
$$
\begin{array}{l}
\displaystyle\partial_{\sigma_j}f^i(x)~:=\sum_{1\leq l\leq r}\sigma^l_{j}(x)
~\partial_{l}f^i(x)=\nabla f^i(x)^{\prime}\sigma_{j}(x)=\sigma_{j}(x)^{\prime}\nabla f^i(x)\\
\\
\Longleftrightarrow  \sigma(x)^{\prime}\nabla f (x)=\nabla_{\sigma}f(x)\Longleftrightarrow
\nabla f (x)^{\prime} \sigma(x)=\nabla_{\sigma}f(x)^{\prime}.
\end{array}
$$
When $\overline{r}=1$, the above formula yields 
$$
\nabla f (x)^{\prime} \sigma_1(x)=\nabla_{\sigma_1}f(x)^{\prime}=\left[\begin{array}{c}
\langle \nabla f^1 (x), \sigma_1(x)\rangle\\ 
\vdots\\
\langle \nabla f^p (x), \sigma_1(x)\rangle
\end{array}\right].
$$
Note that
$$
\nabla_{\sigma_l} \sigma_k=\left(\nabla \sigma_k\right)^{\prime} \sigma_l=
\left[\begin{array}{c}
 (\nabla \sigma^1_k)^{\prime} \sigma_l\\ 
\vdots\\
( \nabla \sigma^r_k)^{\prime} \sigma_l
\end{array}\right]=
\left[\begin{array}{c}
 \nabla_{\sigma_l} \sigma^1_k \\ 
\vdots\\
 \nabla_{\sigma_l} \sigma^r_k
\end{array}\right].
$$
In this notation, the drift function $b_{\sigma}$ of the It\^o formulation  (\ref{diff-def-Ito}) takes the following form
 $$
b_{\sigma}^i:=b^i+\frac{1}{2}~\sum_{1\leq j\leq \overline{r}}~(\nabla \sigma^i_{j})^{\prime}\sigma_{j}=
b^i+\frac{1}{2}~\sum_{1\leq j\leq \overline{r}}~\nabla_{\sigma_j} \sigma^i_{j}.
$$
For any column vector-valued
 smooth function $f=(f^i)_{1\leq i\leq p}$
from $\RR^{r}$ into $\RR^{p}$ we also have the formulae
$$
df(\overline{X}_t)=\nabla_b f (\overline{X}_t)^{\prime}~dt+
\nabla_{\sigma} f (\overline{X}_t)^{\prime} ~d\overline{B}_t
\quad \mbox{\rm and}\quad
df(X_t)=L(f) (X_t)~dt+
\nabla_{\sigma} f (X_t)^{\prime} ~dB_t
$$
with the column vector-valued
 smooth function
$$
L(f)(x)=\nabla_{b_{\sigma}} f (x)^{\prime}+\frac{1}{2}~
  \nabla^2 f(x)^{\prime}~a(x)\quad \mbox{\rm and the matrix}\quad a(x)^{i,j}:=\sum_{1\leq k\leq \overline{r}}~\sigma^i_{k}(x)\sigma^j_{k}(x).
$$
For any time horizon $s\geq 0$ denote by
  $\overline{X}_{s,t}(x)$ and   $X_{s,t}(x)$  the stochastic flows defined for any $t\in [s,\infty[$ and any starting point $\overline{X}_{s,s}(x)=x= X_{s,s}(x)\in \RR^r$ by the equations \eqref{diff-def-over-X} and \eqref{diff-def-Ito}. When $s=0$ we simplify notation and we write   $\overline{X}_{t}(x)$ and   $X_{t}(x)$ instead of   $\overline{X}_{0,t}(x)$ and   $X_{0,t}(x)$.

 For any  $s\geq 0$ denote by
  $Z_{s,t}(x)$  the (deterministic) flow defined for any $t\in [s,\infty[$ and any starting point $Z_{s,s}(x)=x\in\RR^r$ by the ordinary differential equation
$$
dZ_{s,t}(x)=b\left(Z_{s,t}(x)\right)~dt.
$$

\begin{lem}\label{lem-Z-exp}
Assume $(H_b)$. Then for any $x\in\RR^r$ and $s\leq t$ we have the uniform exponential decays
 \begin{equation}\label{tangent-Z}
\Vert\nabla Z_{s,t}(x)\Vert\vee \Vert\nabla^2 Z_{s,t}(x)\Vert\leq c~e^{-\lambda_b(t-s)}
\quad\mbox{and}\quad 
\Vert Z_{s,t}(x)\Vert\leq c~(1+\Vert x\Vert).
\end{equation}
 \end{lem}
 The above lemma is rather well known, but for the convenience of the reader a detailed proof is provided in the Appendix on page~\pageref{lem-Z-exp-proof}. 
 \begin{prop}\label{interpol}
For any $s\leq t$ and any function $f$ on $\RR^r$ we have the interpolation formula
\begin{equation}
\begin{array}{l}
\displaystyle f(\overline{X}_{s,t})=f(Z_{s,t})+\sum_{1\leq j\leq \overline{r}}~\int_s^t~\partial_{\sigma_j}(f\circ Z_{u,t})(Z_{s,u})~d\overline{B}^j_u\\
\\
\hskip3cm+\displaystyle \sum_{1\leq i,j\leq \overline{r}}~\int_{s\leq v\leq u\leq t}
\partial_{\sigma_i}\left(\partial_{\sigma_j}(f\circ Z_{u,t})\circ Z_{v,u}\right)(\overline{X}_{s,v})~d\overline{B}^i_v~d\overline{B}^j_u.
\end{array}\end{equation}
 \end{prop}
  The above interpolation formula is known: see \cite[Theorem 4.1]{dmsumeet} in the context of diffusion flows and we sketch a proof in the Appendix on page~\pageref{interpol-proof}.

 \begin{theo}\label{theo-moments}
 Consider the piecewise smooth anticipative process $  \overline{B}_t$ defined in \eqref{ADT}. In this situation, for any $p\geq 1$ we have the uniform estimate (\ref{unif-intro}).
 \end{theo}
 
 \proof
 To simplify the notation we assume that $r=\overline{r}=1$. Using (\ref{tangent-Z}) we readily check that
 $$
\Vert \left(\partial_{\sigma}Z_{u,t}\right)\Vert\leq c_{\sigma}(1)~e^{-\lambda (t-u)}\quad\mbox{\rm and}\quad
\Vert \partial_{\sigma}\left((\partial_{\sigma}Z_{u,t})\circ Z_{v,u}\right)\Vert\leq c_{\sigma}(2)~e^{-\lambda (t-v)}
$$
with the constants
$$
c_{\sigma}(1):= \Vert\sigma\Vert \quad\mbox{\rm and}\quad c_{\sigma}(2):=\Vert\sigma\Vert~(\Vert \partial\sigma\Vert+\Vert \sigma\Vert).
 $$

For any $t_n\leq t\leq t_{n+1}$ and $1\leq k\leq n$ we set
 $$
\alpha^{(k)}_{s,t}(x):=\frac{1}{\sqrt{t_{k}-t_{k-1}}}~\int_{t_{k-1}}^{s}\left(\partial_{\sigma}Z_{u,t}\right)(Z_{u}(x))~du\quad \mbox{\rm and}\quad
 \goodchi_{t_k}:=\frac{B_{t_{k}}-B_{t_{k-1}}}{\sqrt{t_{k}-t_{k-1}}}.
 $$
Therefore, we have the decomposition
$$
I_t:=\int_0^t\left(\partial_{\sigma}Z_{u,t}\right)(Z_{u})~d\overline{B}_u=I_t^{(1)}+I_t^{(2)}
$$
with the random functions
$$
I^{(1)}_t:=\sum_{1\leq k\leq n}  \goodchi_{t_k}~\alpha^{(k)}_{t_{k},t}\quad \mbox{\rm and}\quad
I^{(2)}_t:= \goodchi_{t_{n+1}}~ \alpha^{(n+1)}_{t,t}.
$$

By Marcinkiewicz-Zygmund's inequality, for any $p\geq 1$ we have
$$
\EE\left(\left(I^{(1)}_t(x)\right)^{2p}\right)^{1/p}\leq c_p~\sum_{1\leq k\leq n}~ (\alpha^{(k)}_{t_{k},t}(x))^2.
$$
 Observe that
 $$
 \begin{array}{l}
  \displaystyle
\sum_{1\leq k\leq n}~ (\alpha^{(k)}_{t_{k},t}(x))^2
=\sum_{1\leq k\leq n}\frac{1}{t_{k}-t_{k-1}}~\left(\int_{t_{k-1}}^{t_k}\left(\partial_{\sigma}Z_{u,t}\right)(Z_{u}(x))~du\right)^2\\
\\
  \displaystyle\leq c~\sum_{1\leq k\leq n}\frac{1}{t_{k}-t_{k-1}}~\left(\int_{t_{k-1}}^{t_k}~e^{-\lambda (t-u)}~du\right)^2\\
  \\
  \displaystyle=(c/\lambda^2)~\sum_{1\leq k\leq n}e^{-2\lambda (t-t_k)}~\frac{1}{t_{k}-t_{k-1}}~\left(1-e^{-\lambda (t_k-t_{k-1})}\right)^2.
   \end{array}
$$
This yields the estimate
$$
\EE\left(\left(I^{(1)}_t(x)\right)^{2p}\right)^{1/p}\leq 
  c_p~\sum_{1\leq k\leq n}e^{-2\lambda (t-t_k)}~(t_k-t_{k-1}).
$$
 We conclude that for any $p\geq 1$
 $$
\sup_{t\geq 0}\sup_{x\in\RR} \left(\EE\left(\left\vert I_t^{(1)}(x)\right\vert^{p}\right)\vee 
 \EE\left(\left\vert I_t^{(2)}(x)\right\vert^{p}\right)\right)<\infty.
 $$
In the same vein, by Proposition~\ref{interpol} for any $t_n\leq t\leq t_{n+1}$  we have the decomposition
$$
\overline{X}_t-Z_t=I_t+J_t$$
with
$$
J_t:=\int_0^t
\int_0^u \partial_{\sigma}\left((\partial_{\sigma}Z_{u,t})\circ Z_{v,u}\right)(\overline{X}_{v})~d\overline{B}_v~d\overline{B}_u=J_t^{(1)}+J^{(2)}_t
$$
and the random functions
 \begin{eqnarray*} 
J_t^{(1)}&:=& \sum_{1\leq k\leq n}\frac{B_{t_{k}}-B_{t_{k-1}}}{t_{k}-t_{k-1}}~\int_{t_{k-1}}^{t_k}
\left(\int_0^u \partial_{\sigma}\left((\partial_{\sigma}Z_{u,t})\circ Z_{v,u}\right)(\overline{X}_{v})~d\overline{B}_v\right)~du\\
\\
J_t^{(2)}&=&\frac{B_{t_{n+1}}-B_{t_{n}}}{t_{n+1}-t_{n}}~\int_{t_{n}}^{t}
\left(\int_0^u \partial_{\sigma}\left((\partial_{\sigma}Z_{u,t})\circ Z_{v,u}\right)(\overline{X}_{v})~d\overline{B}_v\right)~du.
 \end{eqnarray*}
For any $1\leq k\leq n$ and $1\leq l<k$ we set
 \begin{eqnarray*} 
\beta^{(k)}_t(x)&:=&\frac{1}{t_k-t_{k-1}}~\int_{t_{k-1}}^{t_k}
\left(\int_{t_{k-1}}^u \partial_{\sigma}\left((\partial_{\sigma}Z_{u,t})\circ Z_{v,u}\right)(\overline{X}_{v}(x))~dv\right)~du\\
\beta^{(l,k)}_t(x)&:=&\frac{1}{\sqrt{t_k-t_{k-1}}}~\frac{1}{\sqrt{t_l-t_{l-1}}}~\int_{t_{k-1}}^{t_k}
\left(\int_{t_{l-1}}^{t_l} \partial_{\sigma}\left((\partial_{\sigma}Z_{u,t})\circ Z_{v,u}\right)(\overline{X}_{v}(x))~dv\right)~du.
 \end{eqnarray*}
We readily check the almost sure uniform estimates
 \begin{eqnarray*} 
\sup_{x\in\RR} \left\vert\beta^{(k)}_t(x)\right\vert&\leq &\frac{1}{t_k-t_{k-1}}~e^{-\lambda (t-t_k)}~\int_{t_{k-1}}^{t_k}
\int_{t_{k-1}}^u e^{-\lambda (t_k-v)}~dvdu\leq e^{-\lambda (t-t_k)}~(t_k-t_{k-1}) \end{eqnarray*}
 and
  \begin{eqnarray*} 
\sup_{x\in\RR} \left(\beta^{(l,k)}_t(x)\right)^2&\leq &\frac{e^{-2\lambda (t-t_k)}}{t_{k}-t_{k-1}}~  
\frac{e^{-2\lambda (t_k-t_l)}}{t_{l}-t_{l-1}}~\left(\int_{t_{k-1}}^{t_k}
\int_{t_{l-1}}^{t_l} ~e^{-\lambda (t_l-v)}~dv du\right)^2\\
&\leq& e^{-2\lambda (t-t_k)}~(t_{k}-t_{k-1})~  
~e^{-2\lambda (t_k-t_l)} (t_{l}-t_{l-1}).
 \end{eqnarray*}
 
We have the decomposition
 $$
J_t^{(1)}=J^{(1,1)}_t+J^{(1,2)}_t
$$
with
 \begin{eqnarray*} 
J_t^{(1,1)}&:=&\sum_{1\leq k\leq n}~\beta^{(k)}_t~\goodchi_{t_k}^2\\
 J_t^{(1,2)}&:=&\sum_{1\leq k\leq n}
\gamma^{(k-1)}_t~\goodchi_{t_k} 
\quad \mbox{\rm and}\quad\gamma^{(k-1)}_t:=\sum_{1\leq l<k}\beta^{(l,k)}_t~\goodchi_{t_l}.   
 \end{eqnarray*}
One can check that
 $$
 \vert J_t^{(1,1)}(x)\vert\leq \sum_{1\leq k\leq n}~\goodchi_{t_k}^2~e^{-\lambda (t-t_k)}~(t_k-t_{k-1}). 
 $$
 Thus, for any $p\geq 1$ we have
$$
\sup_{x\in\RR}\EE\left( \vert J_t^{(1,1)}(x)\vert^p\right)^{1/p}\leq c_p \sum_{1\leq k\leq n}~e^{-\lambda (t-t_k)}~(t_k-t_{k-1}).
 $$
 Observe that for any $t_{l-1}\leq v\leq t_l$, the random state $
\overline{X}_{v}$ and thus $\beta^{(l,k)}_t$  are measurable w.r.t.~the random vector $\left(
B_{t_{m}}-B_{t_{m-1}}\right)_{1\leq m\leq l}$. Thus the random state
$\gamma_{t_{k-1}}(t)$ is measurable w.r.t.~the random vector $\left(
B_{t_{m}}-B_{t_{m-1}}\right)_{1\leq m<k}$. Applying Burkholder's inequality, for any $p\geq 1$ we have
 \begin{eqnarray*} 
\EE\left(\left( J_t^{(1,2)}(x)\right)^{2p}\right)^{1/p}&\leq &c_p(1)~\EE\left(
\left(\sum_{1\leq k\leq n} \left(\gamma^{(k-1)}_t(x)~\goodchi_{t_k} \right)^2
\right)^p
\right)^{1/p}\\
&\leq &c_p(1)~\sum_{1\leq k\leq n} \EE\left(\left(\gamma^{(k-1)}_t(x)~\goodchi_{t_k} \right)^{2p}\right)^{1/p}
\leq c_p(2)~\sum_{1\leq k\leq n} \EE\left((\gamma^{(k-1)}_t(x))^{2p}\right)^{1/p}.
 \end{eqnarray*}
Similarly, we have
 \begin{eqnarray*} 
\EE\left((\gamma^{(k-1)}_t(x))^{2p}\right)^{1/p}
&\leq &
c_p(3)~\sum_{1\leq l<k} \EE\left(\left(\beta^{(l,k)}_t(x)~\goodchi_{t_l}   
 \right)^{2p}\right)^{1/p}\\
 &\leq &c_p(4)~e^{-2\lambda (t-t_k)}~(t_{k}-t_{k-1})~  \sum_{1\leq l<k} 
~e^{-2\lambda (t_k-t_l)} (t_{l}-t_{l-1})
 \end{eqnarray*} 
 from which we check that
 $$
 \EE\left(\left( J_t^{(1,2)}(x)\right)^{2p}\right)^{1/p}\leq \sum_{1\leq k\leq n}
 e^{-2\lambda (t-t_k)}~(t_{k}-t_{k-1})~  \sum_{1\leq l<k} 
~e^{-2\lambda (t_k-t_l)} (t_{l}-t_{l-1}).
 $$
 We conclude that for any $p\geq 1$
 $$
\sup_{t\geq 0}\sup_{x\in\RR} \left(\EE\left(\left\vert J_t^{(1)}(x)\right\vert^{p}\right)\vee  \EE\left(\left\vert J_t^{(2)}(x)\right\vert^{p}\right)\right)<\infty.
 $$
This ends the proof of the theorem.
 \cqfd

 \section{Homogenous Diffusion Matrices}\label{theo-sigma-ct-proof-i}
 This section is mainly concerned with the proof of Theorem~\ref{theo-sigma-ct}.
We set $\sigma(x)=\Sigma\in \RR^{r\times r}$ and
$$
\Delta_t(x):=\left(X_{t}(x)-\overline{X}_{t}(x)\right)-\Sigma~(B_t-
\overline{B}_t)=\int_0^t\left(b\left(X_{s}(x)\right)-b\left( \overline{X}_{s}(x)\right)\right)~ds.
$$
Thus, we have
\begin{eqnarray*}
\partial_t~ \Vert \Delta_t(x)\Vert^2&=&\partial_t(\Delta_t^{\prime}(x)\Delta_t(x))=2~\Delta_t(x)^{\prime}
\left(b\left(X_{t}(x)\right)-b\left( \overline{X}_{t}(x)\right)\right)\\
&\leq &-2\lambda_b~\Vert X_{t}(x)-\overline{X}_{t}(x)\Vert^2-2
\left(b\left(X_{t}(x)\right)-b\left( \overline{X}_{t}(x)\right)\right)^{\prime}\Sigma(B_t-
\overline{B}_t).
\end{eqnarray*}
Now, we note that
$$
\begin{array}{l}
\displaystyle b(x)-b(y)=\int_0^1~\nabla b(u x+(1-u)y)^{\prime}(x-y)~du\\
\\
\displaystyle\Longrightarrow \Vert b(x)-b(y)\Vert\leq \Vert \nabla b\Vert_2~\Vert x-y\Vert
\quad \mbox{\rm with}\quad
\Vert \nabla b\Vert_2:=\sup_{z\in \RR^r}\Vert \nabla b(z)\Vert_2.
\end{array}
$$
Recalling that $2x^{\prime}y\leq \delta \Vert x\Vert^2+\frac{\Vert y\Vert^2}{\delta}$ for any $\delta>0$ and $x,y\in \RR^r$  we check the estimates
\begin{eqnarray*}
\partial_t~ \Vert \Delta_t(x)\Vert^2&\leq& -2\left(\lambda_b-\delta\Vert \nabla b\Vert_2^2\right)~\Vert X_{t}(x)-\overline{X}_{t}(x)\Vert^2+\frac{1}{2\delta}~\Vert \Sigma\Vert^2_2~\Vert B_t-
\overline{B}_t\Vert^2\\
&=
&-2\left(\lambda_b-\delta\Vert \nabla b\Vert_2^2\right)~\Vert \Delta_t(x)\Vert^2
-4\left(\lambda_b-\delta\Vert \nabla b\Vert_2^2\right)
\Delta_t(x)^{\prime}\Sigma~(B_t-
\overline{B}_t)\\
&&\hskip1cm+\left(\frac{1}{2\delta}
-2\left(\lambda_b-\delta\Vert \nabla b\Vert_2^2\right)\right)~\Vert \Sigma\Vert^2_2~\Vert B_t-
\overline{B}_t\Vert^2.
\end{eqnarray*}
The last assertion follows from the fact that
$$
\begin{array}{l}
\Vert X_{t}(x)-\overline{X}_{t}(x)\Vert^2=\Vert \Delta_t(x)\Vert^2+2~\Delta_t(x)^{\prime}\Sigma~(B_t-
\overline{B}_t)
+\Vert \Sigma\Vert^2_2~\Vert B_t-
\overline{B}_t\Vert.
\end{array}
$$
Choosing $\delta<1\wedge(\lambda_b/\Vert \nabla b\Vert_2^2)$ we find that
\begin{eqnarray*}
\partial_t~ \Vert \Delta_t(x)\Vert^2&\leq
&-2\lambda(\delta)~\Vert \Delta_t(x)\Vert^2+c(\delta)~\Vert B_t-
\overline{B}_t\Vert^2
\end{eqnarray*}
with the parameters
$$
\lambda(\delta):=(1-\delta)~\left(\lambda_b-\delta\Vert \nabla b\Vert_2^2\right)\quad \mbox{\rm and}\quad
c(\delta):=\frac{2}{\delta}~\left(\frac{1}{4}+\lambda(\delta)\right)~\Vert \Sigma\Vert^2_2.
$$
This yields the  uniform estimate
$$
\sup_{x\in\RR^r}\Vert \Delta_t(x)\Vert^2\leq c(\delta)~\int_0^t~e^{-2\lambda(\delta)(t-s)}~
\Vert B_s-
\overline{B}_s\Vert^2
ds.
$$
To conclude, observe that
$$
\sup_{x\in\RR^r}\Vert X_{t}(x)-\overline{X}_{t}(x)\Vert^2 \leq 2\sup_{x\in\RR^r}\Vert \Delta_t(x)\Vert^2+2\Vert \Sigma\Vert_2^2~\Vert B_t-
\overline{B}_t\Vert^2.
$$
\cqfd

\section{Milstein-Type Expansions}
\label{sec:proof_theo_loc}
This section is mainly concerned with the proof of Theorem~\ref{theo-loc-intro}.
\subsection{Local decompositions}
Observe that
$$
\Delta  \overline{X}^i_{t_{n+1}}(x)=
\int_{t_n}^{t_{n+1}} b^i\left( \overline{X}_{t}(x)\right)~dt+\frac{1}{\sqrt{\epsilon}}~\sum_{1\leq k\leq \overline{r}}\int_{t_n}^{t_{n+1}}~ 
\sigma_k^i\left( \overline{X}_{t}(x)\right)~dt~\chi^k_{t_{n+1}}.
$$
Conversely, by integration by parts for any function $f$ we have
$$
\begin{array}{l}
\displaystyle \left(f\left( \overline{X}_{t}(x)\right)-f\left( \overline{X}_{t_n}(x)\right)\right)~dt\\
\\
\displaystyle
=-d\left(\left(f\left( \overline{X}_{t}(x)\right)-f\left( \overline{X}_{t_n}(x)\right)\right)~(t_{n+1}-t)\right)+(t_{n+1}-t)~d\left(f\left( \overline{X}_{t}(x)\right)\right).
\end{array}
$$
This yields the formula
$$
\begin{array}{l}
\displaystyle\int_{t_n}^{t_{n+1}} \left(f\left( \overline{X}_{t}(x)\right)-f\left( \overline{X}_{t_n}(x)\right)\right)~dt\\
\\
=\displaystyle \displaystyle\int_{t_n}^{t_{n+1}} (t_{n+1}-t)~~
\nabla_b f(\overline{X}_{t}(x))~ dt+~\frac{1}{\sqrt{\epsilon}}~\sum_{1\leq k\leq \overline{r}}\int_{t_n}^{t_{n+1}} (t_{n+1}-t)~
\nabla_{\sigma_k} f(\overline{X}_{t}(x))~dt~\chi^k_{t_{n+1}}
\end{array}
$$
from which we check that
$$
\begin{array}{l}
\displaystyle \Delta  \overline{X}^i_{t_{n+1}}(x)=b^i\left( \overline{X}_{t_n}(x)\right)~\epsilon+\sum_{1\leq k\leq \overline{r}}\sigma^i_k\left( \overline{X}_{t_n}(x)\right)\sqrt{\epsilon}~\goodchi_{t_{n+1}}^k~\\
\\
+\displaystyle\int_{t_n}^{t_{n+1}} (t_{n+1}-t)~~
\nabla_b b^i(\overline{X}_{t}(x))~dt+~\frac{1}{\sqrt{\epsilon}}~\sum_{1\leq k\leq \overline{r}}\int_{t_n}^{t_{n+1}} (t_{n+1}-t)~
(\nabla_{\sigma_k} b^i+\nabla_b \sigma^i_k)(\overline{X}_{t}(x))~dt~\chi^k_{t_{n+1}}\\
\\
\displaystyle+~\frac{1}{\epsilon}~\sum_{1\leq k,l\leq \overline{r}}~\int_{t_n}^{t_{n+1}} (t_{n+1}-t)~
\nabla_{\sigma_l} \sigma_k^i(\overline{X}_{t}(x))~dt~~\goodchi_{t_{n+1}}^k~\goodchi_{t_{n+1}}^l.
\end{array}
$$
We can also establish that
$$
\begin{array}{l}
\displaystyle (t_{n+1}-t)~\left(f\left( \overline{X}_{t}(x)\right)-f\left( \overline{X}_{t_n}(x)\right)\right)~dt\\
\\
\displaystyle
=-\frac{1}{2}~d\left(\left(f\left( \overline{X}_{t}(x)\right)-f\left( \overline{X}_{t_n}(x)\right)\right)~(t_{n+1}-t)^2\right)+\frac{1}{2}~(t_{n+1}-t)^2~d\left(f\left( \overline{X}_{t}(x)\right)\right)
\end{array}
$$
and
$$
d\left(f\left( \overline{X}_{t}(x)\right)\right)=\nabla_b f (\overline{X}_{t}(x))~ dt+\frac{1}{\sqrt{\epsilon}}~\sum_{1\leq j\leq \overline{r}}
\nabla_{\sigma_j} f (\overline{X}_{t}(x))~dt~\goodchi_{t_{n+1}}^j.
$$
This yields the formula
$$
\begin{array}{l}
\displaystyle\int_{t_n}^{t_{n+1}}  (t_{n+1}-t)\left(f\left( \overline{X}_{t}(x)\right)-f\left( \overline{X}_{t_n}(x)\right)\right)~dt\\
\\
=\displaystyle \displaystyle\frac{1}{2}\int_{t_n}^{t_{n+1}} (t_{n+1}-t)^2~~
\nabla_b f(\overline{X}_{t}(x)) ~dt+~\frac{1}{2}~\frac{1}{\sqrt{\epsilon}}~\sum_{1\leq j\leq \overline{r}}~\int_{t_n}^{t_{n+1}} (t_{n+1}-t)^2~
\nabla_{\sigma_j} f(\overline{X}_{t}(x))~dt~\goodchi_{t_{n+1}}^j.
\end{array}
$$
Applying the above formula to $f=\nabla_{\sigma_l} \sigma_k^i$ we check that
$$
\begin{array}{l}
\displaystyle \int_{t_n}^{t_{n+1}} (t_{n+1}-t)~
\nabla_{\sigma_l} \sigma_k^i(\overline{X}_{t}(x))~dt
\\
\\
\displaystyle=\frac{1}{2}~\epsilon^2~
\nabla_{\sigma_l} \sigma_k^i(\overline{X}_{t_n}(x))
+\frac{1}{2}~\int_{t_n}^{t_{n+1}} (t_{n+1}-t)^2~
\nabla_b \nabla_{\sigma_l} \sigma_k^i (\overline{X}_{t}(x))~ dt\\
\\
\hskip3cm\displaystyle+\frac{1}{\sqrt{\epsilon}}~\frac{1}{2}~\sum_{1\leq j\leq \overline{r}}\int_{t_n}^{t_{n+1}} (t_{n+1}-t)^2~
\nabla_{\sigma_j} \nabla_{\sigma_l} \sigma_k^i (\overline{X}_{t}(x))~dt~\goodchi_{t_{n+1}}^j.
\end{array}
$$
This yields the decomposition \eqref{delta-overX-intro} with the remainder term
\begin{eqnarray*}
\displaystyle \epsilon^{3/2}~\Delta \overline{R}_{t_{n+1}}  &:=&\int_{t_n}^{t_{n+1}} (t_{n+1}-t)~~
\nabla_b b^i(\overline{X}_{t})~dt\\
&&+~\frac{1}{\sqrt{\epsilon}}~\sum_{1\leq k\leq \overline{r}}\int_{t_n}^{t_{n+1}} (t_{n+1}-t)~
(\nabla_{\sigma_k} b^i+\nabla_b \sigma^i_k)(\overline{X}_{t})~dt~\chi^k_{t_{n+1}}\\
&&
+\frac{1}{2\epsilon}~\sum_{1\leq k,l\leq \overline{r}}\int_{t_n}^{t_{n+1}} (t_{n+1}-t)^2~
\nabla_b \nabla_{\sigma_l} \sigma_k^i (\overline{X}_{t})~ dt~\goodchi_{t_{n+1}}^k~\goodchi_{t_{n+1}}^l\\
&&
\displaystyle+\frac{1}{2\epsilon\sqrt{\epsilon}}~\sum_{1\leq j,k,l\leq \overline{r}}\int_{t_n}^{t_{n+1}} (t_{n+1}-t)^2~
\nabla_{\sigma_j} \nabla_{\sigma_l} \sigma_k^i (\overline{X}_{t})~dt~\goodchi_{t_{n+1}}^j\goodchi_{t_{n+1}}^k~\goodchi_{t_{n+1}}^l.
\end{eqnarray*}
The end of the proof of \eqref{delta-overX-intro} with the r.h.s.~estimate in \eqref{est-r-X} is now a direct consequence of the uniform estimate \eqref{unif-intro}, see also  Theorem~\ref{theo-moments}.

Now we come to the proof of \eqref{delta-X-intro}. We note that
$$
\Delta X_{t_{n+1}}(x)=b_{\sigma}\left(X_{t_n}(x)\right)~\epsilon+\sum_{1\leq k\leq \overline{r}}\sigma_k\left(X_{t_n}(x)\right)\sqrt{\epsilon}~\chi^k_{t_{n+1}}+\Delta \Ra_{t_{n+1}} (x)
$$
with the remainder term
$$
\begin{array}{l}
\displaystyle
\Delta \Ra_{t_{n+1}} (x):=\int_{t_n}^{t_{n+1}}~(b_{\sigma}(X_s(x))-b_{\sigma}(X_{t_n}(x)))~ds+
\int_{t_n}^{t_{n+1}}~(\sigma(X_s(x))-\sigma(X_{t_n}(x)))~dB_s\\
\\
\displaystyle=\int_{t_n}^{t_{n+1}}~(t_{n+1}-t)~L(b_{\sigma})(X_t(x))~dt+
\int_{t_n}^{t_{n+1}}~(t_{n+1}-t)~\nabla_{\sigma}b_{\sigma}(X_t(x))^{\prime}~dB_t\\
\\
\displaystyle+\sum_{1\leq k\leq \overline{r}}\int_{t_n}^{t_{n+1}}~\left(\int_{t_n}^s L(\sigma_k)(X_u(x))~du\right)~dB_s^k+\sum_{1\leq k\leq \overline{r}}
\int_{t_n}^{t_{n+1}}~
\left(\int_{t_n}^s \nabla_{\sigma}\sigma_k(X_u(x))^{\prime}~dB_u\right)~dB_s^k.
\end{array}
$$

%



Observe that
$$
\begin{array}{l}
\displaystyle
\sum_{1\leq k\leq \overline{r}}\int_{t_n}^{t_{n+1}}~
\left(\int_{t_n}^s \nabla_{\sigma}\sigma_k(X_u(x))^{\prime}~dB_u\right)~dB^k_s\\
\\
\displaystyle=\sum_{1\leq k,l\leq \overline{r}}\int_{t_n}^{t_{n+1}}~
\left(\int_{t_n}^s \nabla \sigma_k(X_u(x))^{\prime} \sigma_{l}(X_u(x))~dB_u^l\right)~dB^k_s\\
\\
\displaystyle=\sum_{1\leq k,l\leq \overline{r}} \nabla \sigma_k(X_{t_n}(x))^{\prime} \sigma_{l}(X_{t_n}(x))~\int_{t_n}^{t_{n+1}}~
\int_{t_n}^s~dB_u^l~dB^k_s+\epsilon~\Delta M_{t_{n+1}}(x)
\end{array}$$
with the martingale increment
$$
\Delta M_{t_{n+1}}(x):=\epsilon^{-1}
\sum_{1\leq k,l\leq \overline{r}}
\int_{t_n}^{t_{n+1}}~
\left(\int_{t_n}^s \left(\nabla_{\sigma_l}\sigma_k(X_{t_n,u}(X_{t_n}))-\nabla_{\sigma_l}\sigma_k(X_{t_n})\right)~dB^l_u\right)~dB^k_s.
$$
Using the commutation property
$$
\nabla \sigma_k(x)^{\prime} \sigma_{l}(x)=
\nabla \sigma_l(x)^{\prime} \sigma_{k}(x)
$$
we check that
$$
\begin{array}{l}
\displaystyle
\sum_{1\leq k\leq \overline{r}}\int_{t_n}^{t_{n+1}}~
\left(\int_{t_n}^s \nabla_{\sigma}\sigma_k(X_u(x))^{\prime}~dB_u\right)~dB^k_s\\
\\
\displaystyle=\frac{1}{2}~\epsilon~\sum_{1\leq k,l\leq \overline{r}}
\nabla_{\sigma_l}\sigma_k(X_{t_n}(x))\left(\goodchi_{t_{n+1}}^l\goodchi_{t_{n+1}}^k-1_{k=l}\right)+\epsilon~\Delta M_{t_{n+1}}(x).
\end{array}$$
The last assertion follows as
$$
\begin{array}{l}
\displaystyle
\int_{t_n}^{t_{n+1}}~(B^l_s-B^l_{t_n})~dB^k_s+\int_{t_n}^{t_{n+1}}~(B^k_s-B^k_{t_n})~dB^l_s\\
\\
\displaystyle=\left((B^l_{t_{n+1}}-B^l_{t_n})
(B^k_{t_{n+1}}-B^k_{t_n})-1_{k=l}(t_{n+1}-t_n)\right)=
\epsilon~\left(\goodchi_{t_{n+1}}^l\goodchi_{t_{n+1}}^k-1_{k=l}\right).\end{array}$$
This ends the proof of \eqref{delta-X-intro} with the remainder term
\begin{eqnarray*}
\epsilon^{3/2}~\Delta R_{t_{n+1}}
&:=&\int_{t_n}^{t_{n+1}}~(t_{n+1}-t)~L(b_{\sigma})(X_t)~dt+
\int_{t_n}^{t_{n+1}}~(t_{n+1}-t)~\nabla_{\sigma}b_{\sigma}(X_t)^{\prime}~dB_t\\
&&+\sum_{1\leq k\leq \overline{r}}\int_{t_n}^{t_{n+1}}~\left(\int_{t_n}^s L(\sigma_k)(X_u)~du\right)~dB_s^k.
\end{eqnarray*}
Using (\ref{hyp-mom}), for any $p\geq 1$ we check the uniform estimate
$$
\sup_x\EE\left(\Vert \int_{t_n}^{t_{n+1}}~(t_{n+1}-t)~L(b_{\sigma})(X_t(x))~dt\Vert^p\right)^{1/p}\leq c_p(1)~\epsilon^2~(1+\vert x\vert)
$$
as well as
$$
\begin{array}{l}
\displaystyle
\EE\left(\Vert\int_{t_n}^{t_{n+1}}~(t_{n+1}-t)~\nabla_{\sigma}b_{\sigma}(X_t(x))^{\prime}~dB_t\Vert^{2p}\right)^{1/p}\leq c_{p}(2)~\epsilon^{3}~
\end{array}$$
and
$$
\begin{array}{l}
\displaystyle
\EE\left(\Vert\int_{t_n}^{t_{n+1}}~\left(\int_{t_n}^t L(\sigma)(X_u(x))~du\right)~dB_t\Vert^{2p}\right)^{1/p}\leq 
c_{p}(3)~\epsilon^3~(1+\vert x\vert)^2.
\end{array}
$$
For any $1\leq j,k,l\leq \overline{r}$ we also have
$$
\begin{array}{l}
\displaystyle
\EE\left(\chi^j_{t_{n+1}}
\int_{t_n}^{t_{n+1}}~
\left(\int_{t_n}^s \left(\nabla_{\sigma_l}\sigma_k(X_{u})-\nabla_{\sigma_l}\sigma_k(X_{t_n})\right)~dB^l_u\right)~dB^k_s~
\vert~\Fa_{t_n}\right)\\
\\
=\displaystyle 1_{j=k}
\int_{t_n}^{t_{n+1}}~\EE\left(
\left(\int_{t_n}^s \left(\nabla_{\sigma_l}\sigma_k(X_{u})-\nabla_{\sigma_l}\sigma_k(X_{t_n})\right)~dB^l_u\right)~\vert~\Fa_{t_n}\right)ds=0
\end{array}
$$
yielding proof of \eqref{est-r-M}. For any $1\leq i\leq r$ we also have
$$
\begin{array}{l}
\displaystyle
\EE\left(\left(
\int_{t_n}^{t_{n+1}}~
\left(\int_{t_n}^s \left(\nabla_{\sigma_l}\sigma_k^i(X_{u})-\nabla_{\sigma_l}\sigma_k^i(X_{t_n})\right)~dB^l_u\right)~dB^k_s\right)^{2p}~|~\Fa_{t_n}\right)^{1/p}\\
\\
\leq \displaystyle c_{p}(1)~
\EE\left(\left(\int_{t_n}^{t_{n+1}}~
\left(\int_{t_n}^s \left(\nabla_{\sigma_l}\sigma_k^i(X_{u})-\nabla_{\sigma_l}\sigma_k^i(X_{t_n})\right)~dB^l_u\right)^2~ds\right)^{p}~|~\Fa_{t_n}\right)^{1/p}
\end{array}
$$
$$
\begin{array}{l}
\displaystyle
\leq \displaystyle c_{p}(1)~\int_{t_n}^{t_{n+1}}~
\EE
\left(\left(\int_{t_n}^s \left(\nabla_{\sigma_l}\sigma_k^i(X_{u})-\nabla_{\sigma_l}\sigma_k^i(X_{t_n})\right)~dB^l_u\right)^{2p}~|~\Fa_{t_n}\right)^{1/p}~ds\\
\\
 \displaystyle \leq c_{p}(2)~\int_{t_n}^{t_{n+1}}~
 \EE\left(
\left(\int_{t_n}^s \left(\nabla_{\sigma_l}\sigma_k^i(X_{u})-\nabla_{\sigma_l}\sigma_k^i(X_{t_n})\right)^2 ds\right)^{p}~|~\Fa_{t_n}\right)^{1/p}~ds
\leq c_{p}(3)~\epsilon^2
\end{array}
$$
which completes the proof of the estimate \eqref{est-r-X}.
\cqfd

\subsection{Time Uniform Estimates}\label{theo-unif-2-proof}
This section is dedicated to the proof of Theorem~\ref{theo-unif-2}.
We start by noting that
$$
\begin{array}{l}
\displaystyle\Delta X_{t_{n}}-\Delta \overline{X}_{t_{n}} \\
\\
\displaystyle=\left(b_{\sigma}\left(X_{t_{n-1}}\right)-b_{\sigma}\left(\overline{X}_{t_{n-1}}\right)\right)~\epsilon+\left(\sigma\left(X_{t_{n-1}}\right)-\sigma\left(\overline{X}_{t_{n-1}}\right)\right)\sqrt{\epsilon}~\goodchi_{t_{n}}+\epsilon~\Delta M_{t_{n}}\\
\\\displaystyle+2^{-1}~\sum_{1\leq k,l\leq \overline{r}}
\left(\nabla_{\sigma_l}\sigma_k(X_{t_{n-1}})-\nabla_{\sigma_l}\sigma_k(\overline{X}_{t_{n-1}})\right)\left(\goodchi_{t_{n}}^l\goodchi_{t_{n}}^k-1_{k=l}\right)~\epsilon+\epsilon^{3/2}~\left(\Delta R_{t_{n}} -\Delta \overline{R}_{t_{n}} \right).
\end{array}
$$
This yields the formula
$$
\begin{array}{l}
\displaystyle\EE\left(\left(X_{t_{n-1}}-\overline{X}_{t_{n-1}}\right)^{\prime}\left(\Delta X_{t_{n}}-\Delta \overline{X}_{t_{n}}\right)\right) \\
\\
\displaystyle=\EE\left(\left(X_{t_{n-1}}-\overline{X}_{t_{n-1}}\right)^{\prime}\left(b_{\sigma}\left(X_{t_{n-1}}\right)-b_{\sigma}\left(\overline{X}_{t_{n-1}}\right)\right)\right)~\epsilon+\epsilon^{3/2}~\EE\left(\left(X_{t_{n-1}}-\overline{X}_{t_{n-1}}\right)^{\prime}\left(\Delta R_{t_{n}} -\Delta \overline{R}_{t_{n}} \right)\right).
\end{array}
$$
In addition, after some elementary but tedious calculations we have the decomposition
$$
\begin{array}{l}
\displaystyle\Vert \Delta X_{t_{n}}-\Delta \overline{X}_{t_{n}}\Vert^2-\Vert\left(\sigma\left(X_{t_{n-1}}\right)-\sigma\left(\overline{X}_{t_{n-1}}\right)\right)\goodchi_{t_{n}}\Vert^2~\epsilon \\
\\
\displaystyle=N_{t_n}+\epsilon^2~\Vert\Delta M_{t_{n}}\Vert^2+\sum_{k=3}^6\epsilon^{k/2}~U_{t_n}^{(k)}\end{array}
$$
with the centered random function
\begin{eqnarray*}
N_{t_n}&:=&2\epsilon^{3/2}~(\Delta M_{t_{n}})^{\prime}\left(\sigma\left(X_{t_{n-1}}\right)-\sigma\left(\overline{X}_{t_{n-1}}\right)\right)~\goodchi_{t_{n}}+2\epsilon^2~\left(b_{\sigma}\left(X_{t_{n-1}}\right)-b_{\sigma}\left(\overline{X}_{t_{n-1}}\right)\right)^{\prime}\Delta M_{t_{n}}\\
&&+\epsilon^2~\sum_{1\leq k,l\leq \overline{r}}~(\Delta M_{t_{n}})^{\prime}~
\left(\nabla_{\sigma_l}\sigma_k(X_{t_{n-1}})-\nabla_{\sigma_l}\sigma_k(\overline{X}_{t_{n-1}})\right)\left(\goodchi_{t_{n}}^l\goodchi_{t_{n}}^k-1_{k=l}\right)
\end{eqnarray*}
and the collection of random functions $U_{t_n}^{(p)}$ defined by
\begin{eqnarray*}
U_{t_n}^{(3)}&:=&2~\left(b_{\sigma}\left(X_{t_{n-1}}\right)-b_{\sigma}\left(\overline{X}_{t_{n-1}}\right)\right)^{\prime}\left(\sigma\left(X_{t_{n-1}}\right)-\sigma\left(\overline{X}_{t_{n-1}}\right)\right)\goodchi_{t_{n}}\\
U_{t_n}^{(4)}&:=&\Vert b_{\sigma}\left(X_{t_{n-1}}\right)-b_{\sigma}\left(\overline{X}_{t_{n-1}}\right)\Vert^2+2\left(\Delta R_{t_{n}} -\Delta \overline{R}_{t_{n}} \right)^{\prime}
\left(\sigma\left(X_{t_{n-1}}\right)-\sigma\left(\overline{X}_{t_{n-1}}\right)\right)~\goodchi_{t_{n}}\\
&&+2^{-2}~\Vert\sum_{1\leq k,l\leq \overline{r}}
\left(\nabla_{\sigma_l}\sigma_k(X_{t_{n-1}})-\nabla_{\sigma_l}\sigma_k(\overline{X}_{t_{n-1}})\right)\left(\goodchi_{t_{n}}^l\goodchi_{t_{n}}^k-1_{k=l}\right)\Vert^2~\\
&&
+2\sum_{1\leq k,l\leq \overline{r}}\left(b_{\sigma}\left(X_{t_{n-1}}\right)-b_{\sigma}\left(\overline{X}_{t_{n-1}}\right)\right)^{\prime}
\left(\nabla_{\sigma_l}\sigma_k(X_{t_{n-1}})-\nabla_{\sigma_l}\sigma_k(\overline{X}_{t_{n-1}})\right)\left(\goodchi_{t_{n}}^l\goodchi_{t_{n}}^k-1_{k=l}\right)
\end{eqnarray*}
and
\begin{eqnarray*}
U^{(5)}_{t_n}&:=&2~\left(b_{\sigma}\left(X_{t_{n-1}}\right)-b_{\sigma}\left(\overline{X}_{t_{n-1}}\right)\right)^{\prime}\left(\Delta R_{t_{n}} -\Delta \overline{R}_{t_{n}} \right)+2~(\Delta M_{t_{n}})^{\prime}~\left(\Delta R_{t_{n}} -\Delta \overline{R}_{t_{n}} \right)\\
&&+\sum_{1\leq k,l\leq \overline{r}}~\left(\Delta R_{t_{n}} -\Delta \overline{R}_{t_{n}} \right)^{\prime}\left(\nabla_{\sigma_l}\sigma_k(X_{t_{n-1}})-\nabla_{\sigma_l}\sigma_k(\overline{X}_{t_{n-1}})\right)\left(\goodchi_{t_{n}}^l\goodchi_{t_{n}}^k-1_{k=l}\right)
\\
U^{(6)}_{t_n}&:=&\Vert\Delta R_{t_{n}} -\Delta \overline{R}_{t_{n}} \Vert^2.
\end{eqnarray*}
Combining \eqref{unif-intro}
and \eqref{hyp-mom} with  \eqref{est-r-X} one has that
$$
\EE\left(\vert U_{t_n}^{(k)}(x)\vert^p\right)^{1/p}\leq c_{p}~(1+\Vert x\Vert^2).
$$
This clearly implies that
$$
\left\vert
\EE\left(\Vert \Delta X_{t_{n}}(x)-\Delta \overline{X}_{t_{n}}(x)\Vert^2\right)-\EE\left(\Vert\sigma\left(X_{t_{n-1}}\right)-\sigma\left(\overline{X}_{t_{n-1}}\right)\Vert^2\right)~\epsilon\right\vert \leq c~\epsilon^{3/2}~(1+\Vert x\Vert^2).
$$
In order to conclude, we have
$$
\begin{array}{l}
\displaystyle \Vert X_{t_{n}}-\overline{X}_{t_{n}}\Vert^2-\Vert  X_{t_{n-1}}- \overline{X}_{t_{n-1}}\Vert^2\\
\\
=\displaystyle 2\left(X_{t_{n-1}}-\overline{X}_{t_{n-1}}\right)^{\prime} \left(\Delta X_{t_{n}}-\Delta \overline{X}_{t_{n}}\right)+\Vert \Delta X_{t_{n}}-\Delta \overline{X}_{t_{n}}\Vert^2
\end{array}
$$
and set
$$
I_n(x):=\EE\left(\Vert X_{t_{n}}(x)-\overline{X}_{t_{n}}(x)\Vert^2\right).
$$
Using the fact that
$$
\begin{array}{l}
\displaystyle\EE\left(\left(X_{t_{n-1}}-\overline{X}_{t_{n-1}}\right)^{\prime}\left(\Delta X_{t_{n}}-\Delta \overline{X}_{t_{n}}\right)\right) \\
\\
\displaystyle=\EE\left(\left(X_{t_{n-1}}-\overline{X}_{t_{n-1}}\right)^{\prime}\left(b_{\sigma}\left(X_{t_{n-1}}\right)-b_{\sigma}\left(\overline{X}_{t_{n-1}}\right)\right)\right)~\epsilon+\epsilon^{3/2}~\EE\left(\left(X_{t_{n-1}}-\overline{X}_{t_{n-1}}\right)^{\prime}\left(\Delta R_{t_{n}} -\Delta \overline{R}_{t_{n}} \right)\right)
\end{array}
$$
we check the estimate
$$
\begin{array}{l}
\displaystyle I_n(x)-I_{n-1}(x)
\leq 2~\EE\left(\left(X_{t_{n-1}}(x)-\overline{X}_{t_{n-1}}(x)\right)^{\prime}\left(b_{\sigma}\left(X_{t_{n-1}}(x)\right)-b_{\sigma}\left(\overline{X}_{t_{n-1}}(x)\right)\right)\right)~\epsilon\\
\\
 \displaystyle\hskip3cm+\EE\left(\Vert\sigma\left(X_{t_{n-1}}(x)\right)-\sigma\left(\overline{X}_{t_{n-1}}(x)\right)\Vert^2\right)~\epsilon\\
\\
 \displaystyle\hskip3cm+\epsilon^{3/2}~\EE\left(\left(X_{t_{n-1}}(x)-\overline{X}_{t_{n-1}}(x)\right)^{\prime}\left(\Delta R_{t_{n}} (x)-\Delta \overline{R}_{t_{n}}(x) \right)\right)+ c~\epsilon^{3/2}~(1+\Vert x\Vert^2).
\end{array}
$$
Whenever  $(H_{\sigma})$ is satisfied, we have
$$
\begin{array}{l}
\displaystyle I_n(x)-I_{n-1}(x)
\leq \displaystyle-2~\lambda_{\sigma}~I_{n-1}(x)~\epsilon\\
\\
\hskip3cm+\epsilon^{3/2}~\EE\left(\left(X_{t_{n-1}}(x)-\overline{X}_{t_{n-1}}(x)\right)^{\prime}\left(\Delta R_{t_{n}} (x)-\Delta \overline{R}_{t_{n}}(x) \right)\right)+ c_1~\epsilon^{3/2}~(1+\Vert x\Vert^2).
\end{array}
$$
Recalling that $2xy\leq \delta x^2+{y^2}/{\delta}$ for any $\delta>0$ and $x,y\in \RR$  and choosing $0<\delta<\lambda_{\sigma}$ we have
$$
\begin{array}{l}
\displaystyle
2\epsilon^{1/2}~\EE\left(\left(X_{t_{n-1}}(x)-\overline{X}_{t_{n-1}}(x)\right)^{\prime}\left(\Delta R_{t_{n}}(x) -\Delta \overline{R}_{t_{n}}(x) \right)\right)\\
\\
\displaystyle\leq 2\delta
\EE\left(\Vert X_{t_{n-1}}-\overline{X}_{t_{n-1}}\Vert^2\right)+c_2\epsilon (1+\Vert x\Vert^2)/\delta.
\end{array}
$$
As a result one has the estimate
$$
\begin{array}{l}
\displaystyle  I_n(x)-I_{n-1}(x)
\leq \displaystyle-2~(\lambda_{\sigma}-\delta)~I_{n-1}~\epsilon
+(c_1~\epsilon^{1/2}+c_2\epsilon/\delta)~(1+\Vert x\Vert^2)
~\epsilon.
\end{array}
$$
We conclude that  for any $0<\delta<1$ we have
$$
\begin{array}{l}
\displaystyle  I_n(x)
\leq \displaystyle(1-2~\lambda_{\sigma}(\delta)\epsilon)~I_{n-1}
+c_{\delta}(\epsilon)~\epsilon^{1/2}~~(1+\Vert x\Vert^2)
~\epsilon
\end{array}
$$
with
$$
\lambda_{\sigma}(\delta):=
(1-\delta/\lambda_{\sigma})\lambda_{\sigma}\quad \mbox{\rm and}\quad
c_{\delta}(\epsilon):=c~(1+\epsilon^{1/2}/\delta).
$$
Therefore, it follows that
$$
I_n(x)\leq c_{\delta}(\epsilon)~\epsilon^{1/2}~(1+\Vert x\Vert)^2\sum_{0\leq k<n}(1-2~\lambda_{\sigma}(\delta)\epsilon)^k~\epsilon
$$
and therefore
$$
\sup_{n\geq 0}I_n(x)\leq ~\epsilon^{1/2}~(1+\Vert x\Vert)^2~\frac{c_{\delta}(\epsilon)}{2~\lambda_{\sigma}(\delta)}
$$
and the proof of \eqref{unif-X-n} is completed.
\cqfd

\section{Wong-Zakai Integration}\label{theo-unif-int-strato-proof-1}

This section is mainly concerned with the proof of Theorem~\ref{theo-unif-int-strato}.
 To simplify the notation we assume that $r=\overline{r}=1$.

\subsection{Local decompositions} 


Consider the stochastic integrals indexed by $t_n\leq u\leq t_{n+1}\leq t$ and defined by
\begin{eqnarray*}
 \Ia_{t_n,u}(F)(x)&:=&\int_{t_n}^{u}~\left(F_{t_n,t}(x)-F_{s,t}(X_{t_n,s}(x))\right)~dB_s\\
 \Ja_{t_n,u}(F)(x)&:=&\int_{t_n}^{u}~\frac{u-s}{t_{n+1}-t_n}~
\left(\sigma
\partial F_{s,t}\right)(X_{t_n,s}(x))~dB_s\\
 \Ka_{t_n,u}(F)(x)&:=&\int_{t_n}^{u}~\frac{u-s}{t_{n+1}-t_n}~
\left(\sigma
\partial F_{s,t}\right)(X_{t_n,s}(x))~ds\\
\La_{t_n,u}(F)(x)&:=&\int_{t_n}^{u}~\frac{u-s}{t_{n+1}-t_n}~
((\partial_s+L)
F_{s,t})(X_{t_n,s}(x))~ds.
\end{eqnarray*}
We also consider the collection of integrals
\begin{eqnarray*}
\Ca_{t_n,u}(F)(x)&:=&
\int_{t_n}^{u}~\frac{u-s}{t_{n+1}-t_n}~
\left(\left(\sigma
\partial F_{s,t}\right)(X_{t_n,s}(x))-\left(\sigma
\partial F_{t_n,t}\right)(x)\right)ds\\
\Da_{t_n,u}(F)(x)&:=&\frac{1}{2}
\int_{t_n}^{u}~\left(\left(\sigma
\partial F_{t_n,t}\right)(x)-\left(\sigma\partial F_{s,t}\right)(X_{t_n,s}(x))\right)~ds
\end{eqnarray*}
and we set
\begin{eqnarray*}
 \Ma_{t_n,u}(F)(x)&:=&\left(B_{u}-B_{t_n}\right)~ \Ja_{t_n,u}(F)(x)- \Ka_{t_n,u}(F)(x)\\
 \Na_{t_n,u}(F)(x)&:=&\left(B_u-B_{t_n}\right)~\La_{t_n,u}(F)(x)+\Ca_{t_n,u}(F)(x).
\end{eqnarray*}

\begin{prop}\label{decomp-1}
For any $t_n\leq u\leq t_{n+1}\leq t$ we  we have the decomposition
\begin{equation}\label{int-over-B}
\begin{array}{l}
\displaystyle\int_{t_n}^{u}~F_{s,t}(X_{t_n,s}(x))~d\overline{B}_s\\
\\
\displaystyle=\left(F_{t_n,t}(x)~\left(B_u-B_{t_n}\right)~+\frac{1}{2}\left(\sigma
\partial F_{t_n,t}\right)(x)~(u-t_n)\right)~\frac{u-t_n}{t_{n+1}-t_n}+ \overline{R}_{t_n,u}(F)(x).
\end{array}
\end{equation}
In the above display, $\overline{R}_{t_n,u}(F)(x)$ stands for the remainder term
$$
\overline{R}_{t_n,u}(F)(x)= \Ma_{t_n,u}(F)(x)+ \Na_{t_n,u}(F)(x)+\overline{\Oa}_{t_n,u}(F)(x)
$$
with
\begin{eqnarray*}
 \overline{\Oa}_{t_n,u}(F)(x)&:=&\left(B_{t_{n+1}}-B_{u}\right)~\left(\La_{t_n,u}(F)(x)+\Ja_{t_n,u}(F)(x)+F_{t_n,t}(x)~\frac{u-t_n}{t_{n+1}-t_n}\right).
\end{eqnarray*}
\end{prop}
\proof
Observe that
\begin{eqnarray*}
\int_{t_n}^{u}~F_{t_n,t}(x)~d\overline{B}_s&=&F_{t_n,t}(x)~\frac{u-t_n}{t_{n+1}-t_n}~\left((B_{t_{n+1}}-B_{u})+(B_u-B_{t_n})\right)\\
\displaystyle\int_{t_n}^{u}~\left(F_{s,t}(X_{t_n,s}(x))-F_{t_n,t}(x)\right)~d\overline{B}_s&=&\int_{t_n}^{u}~\left(F_{s,t}(X_{t_n,s}(x))-F_{t_n,t}(x)\right)~\frac{ds}{t_{n+1}-t_n}~\left(B_{t_{n+1}}-B_{t_n}\right).
\end{eqnarray*}
By integration by parts, we have
$$
\begin{array}{l}
\displaystyle
\left(F_{s,t}(X_{t_n,s}(x))-F_{t_n,t}(x)\right)~\frac{ds}{t_{n+1}-t_n}\\
\\
\displaystyle=-d_s\left(\left(F_{s,t}(X_{t_n,s}(x))-F_{t_n,t}(x)\right)~\frac{u-s}{t_{n+1}-t_n}\right)
+\frac{u-s}{t_{n+1}-t_n}~d_s\left(F_{s,t}(X_{t_n,s}(x))\right)
\end{array}$$
with
$$
\begin{array}{l}
\displaystyle
d_s\left(F_{s,t}(X_{t_n,s}(x))\right)\\
\\
\displaystyle=\left(\left(\partial_s
F_{s,t}\right)(X_{t_n,s}(x))+L\left(F_{s,t}\right)(X_{t_n,s}(x))\right)~ds+\left(\sigma
\partial F_{s,t}\right)(X_{t_n,s}(x))~dB_s.
\end{array}$$
This implies that
$$
\begin{array}{l}
\displaystyle \int_{t_n}^{u}~\left(F_{s,t}(X_{t_n,s}(x))-F_{t_n,t}(x)\right)~d\overline{B}_s\\
\\
\displaystyle =\La_{t_n,u}(F)(x)~\left(B_{t_{n+1}}-B_{t_n}\right) + \Ja_{t_n,u}(F)(x)~\left(B_{t_{n+1}}-B_{t_n}\right)\\
\\
\displaystyle =\left(B_{t_{n+1}}-B_{u}\right) \left(\La_{t_n,u}(F)(x)+ \Ja_{t_n,u}(F)(x)\right) +\Na_{t_n,u}(F)(x)+ \Ma_{t_n,u}(F)(x)\\
\\
\displaystyle\hskip3cm+\Ka_{t_n,u}(F)(x)-   \Ca_{t_n,u}(F)(x).
 \end{array}
$$
Finally note that
$$
\displaystyle\Ka_{t_n,u}(F)(x)-   \Ca_{t_n,u}(F)(x)
=\frac{1}{2}~\left(\sigma
\partial F_{t_n,t}\right)(x)~\frac{(u-t_n)^2}{t_{n+1}-t_n}
$$
from which we conclude.
\cqfd

\begin{prop}\label{prop-int-over-B-2}
For any $t_n\leq u \leq t_{n+1}\leq t$  we have the
\begin{equation}\label{int-over-B-2}
\begin{array}{l}
\displaystyle \int_{t_n}^{u}~F_{s,t}(X_{t_n,s}(x))~d\overline{B}_s=
\int_{t_n}^{u}~F_{s,t}(X_{t_n,s}(x))~dB_s+\frac{1}{2}
\int_{t_n}^{u}~\left(\sigma\partial F_{s,t}\right)(X_{t_n,s}(x))~ds+R_{t_n,u}(F)(x).
\end{array}
\end{equation}
In the above display, ${R}_{t_n,u}(F)(x)$ stands for the remainder term
$$
\begin{array}{l}
\displaystyle R_{t_n,u}(F)(x):=\Ma_{t_n,u}(F)(x)+ \Na_{t_n,u}(F)(x)+ \Ia_{t_n,u}(F)(x)+ \Da_{t_n,u}(F)(x)+{\Oa}_{t_n,u}(F)(x)\end{array}
$$
with
$$
\begin{array}{l}
\displaystyle
{\Oa}_{t_n,u}(F)(x):=-\frac{1}{2}\left(\sigma
\partial F_{t_n,t}\right)(x)~(u-t_n)~\frac{t_{n+1}-u}{t_{n+1}-t_n}+\left(\La_{t_n,u}(F)(x)+\Ja_{t_n,u}(F)(x)\right)\left(B_{t_{n+1}}-B_{u}\right)\\
\\
\hskip3cm\displaystyle+F_{t_n,t}(x)~\left(\frac{u-t_n}{t_{n+1}-t_n}\left(B_{t_{n+1}}-B_{u}\right)~
-\left(B_u-B_{t_n}\right)~\frac{t_{n+1}-u}{t_{n+1}-t_n}\right).
 \end{array}
$$
\end{prop}
\proof
Applying Proposition~\ref{decomp-1} we have
$$
\begin{array}{l}
\displaystyle\int_{t_n}^{u}~F_{s,t}(X_{t_n,s}(x))~d\overline{B}_s=\int_{t_n}^{u}~F_{s,t}(X_{t_n,s}(x))~dB_s+\frac{1}{2}
\int_{t_n}^{u}~\left(\sigma\partial F_{s,t}\right)(X_{t_n,s}(x))~ds\\
\\
\hskip4cm\displaystyle+ \overline{R}_{t_n,u}(F)(x)
+\Ia_{t_n,u}(F)(x)+\Da_{t_n,u}(F)(x)\\
\\
\hskip4cm\displaystyle+
\left(F_{t_n,t}(x)~\left(B_u-B_{t_n}\right)~+\frac{1}{2}\left(\sigma
\partial F_{t_n,t}\right)(x)~(u-t_n)\right)~\left(\frac{u-t_n}{t_{n+1}-t_n}-1\right)\end{array}.
$$
Finally we have
$$
\begin{array}{l}
\displaystyle
{\Oa}_{t_n,u}(F)(x)= \overline{\Oa}_{t_n,u}(F)(x)+\left(F_{t_n,t}(x)~\left(B_u-B_{t_n}\right)~+\frac{1}{2}\left(\sigma
\partial F_{t_n,t}\right)(x)~(u-t_n)\right)~\frac{u-t_{n+1}}{t_{n+1}-t_n}\\
\\
\displaystyle=\left(\La_{t_n,u}(F)(x)+\Ja_{t_n,u}(F)(x)\right)\left(B_{t_{n+1}}-B_{u}\right)-\frac{1}{2}\left(\sigma
\partial F_{t_n,t}\right)(x)~(u-t_n)~\frac{t_{n+1}-u}{t_{n+1}-t_n}\\
\\
\displaystyle+F_{t_n,t}(x)~\left(\left(B_{t_{n+1}}-B_{u}\right)~\frac{u-t_n}{t_{n+1}-t_n}-\left(B_u-B_{t_n}\right)~\frac{t_{n+1}-u}{t_{n+1}-t_n}\right).
 \end{array}
$$
\cqfd

We remark that for any $t_n\leq u \leq t_{n+1}\leq t$ we have
$$
{\Oa}_{t_n,t_{n+1}}(F)(x)=0\quad \mbox{\rm and}\quad
\EE\left({\Oa}_{t_n,u}(F)(x)\right)=-\frac{1}{2}\left(\sigma
\partial F_{t_n,t}\right)(x)~(u-t_n)~\frac{t_{n+1}-u}{t_{n+1}-t_n}
$$
and
$$
\begin{array}{l}
\displaystyle \EE\left(R_{t_n,u}(F)(x)\right)=\EE\left( \left(B_u-B_{t_n}\right)~\La_{t_n,u}(F)(x)\right)+\EE\left(\Ca_{t_n,u}(F)(x)+\Da_{t_n,u}(F)(x)\right)\\
\\
\displaystyle\hskip3cm-\frac{1}{2}\left(\sigma
\partial F_{t_n,t}\right)(x)~(u-t_n)~\frac{t_{n+1}-u}{t_{n+1}-t_n}.\end{array}
$$

\subsection{Some Quantitative Estimates}

\begin{lem}\label{lem-I-int}
For any $t_n\leq u\leq t_{n+1}\leq t$ and $p\geq 1$ we have
$$
 \EE\left(\vert\Ia_{t_n,u}(F)(x)\vert ^{p}\right)^{1/p}\\
\\
\displaystyle\leq  c_p~e^{-\lambda(t-u)}~\left(1+\vert x\vert~\epsilon^{1/2}~\right)~(u-t_n)
$$
and
$$
 \EE\left(\vert \Da_{t_n,u}(F)(x)\vert^{p}\right)^{1/p}
\leq c_p~\epsilon^{1/2}~e^{-\lambda(t-u)}~\left(1+
 \vert x\vert~\epsilon^{1/2}~\right)~(u-t_n).
$$

\end{lem}
\proof
Our regularity conditions ensure that
 for any $u\leq s\leq t$ we have
\begin{equation}\label{contiuity-cond}
\begin{array}{l}
\displaystyle \EE\left(\vert F_{u,t}(x)-F_{s,t}(X_{u,s}(x))\vert^{p}\right)^{1/p}\vee \EE\left(\vert (\sigma\partial F_{u,t})(x)-(\sigma\partial F_{s,t})(X_{u,s}(x))\vert^{p}\right)^{1/p}~\\
\\
\displaystyle \leq  c_p~e^{-\lambda(t-s)}~\left(
 (1+\vert x\vert)~(s-u)+\sqrt{s-u}\right).
\end{array}
\end{equation}
This implies that
$$
\begin{array}{l}
\displaystyle \EE\left(\left(\Ia_{t_n,u}(F)(x)\right)^{2p}\right)^{1/p}\\
\\
\displaystyle\leq c_p(1)~\int_{t_n}^{u}~\EE\left(\left(F_{t_n,t}(x)-F_{s,t}(X_{t_n,s}(x))\right)^{2p}\right)^{1/p}~ds\\
\\
\displaystyle\leq c_p(2)~e^{-2\lambda(t-u)}~\int_{t_n}^{u}~\left(\left(
 (1+\vert x\vert^2)~(s-t_n)^2+(s-t_n)\right)\right)~ds.
\end{array}
$$
This yields the estimate
$$
\begin{array}{l}
\displaystyle \EE\left(\left(\Ia_{t_n,u}(F)(x)\right)^{2p}\right)^{1/(2p)}\leq  c_p(3)~e^{-\lambda(t-u)}~\left(1+
 (1+\vert x\vert)~(u-t_n)^{1/2}\right)~(u-t_n).
\end{array}
$$
In the same vein, we have
$$
\begin{array}{l}
\displaystyle \EE\left(\vert\Da_{t_n,u}(F)(x)\vert^{p}\right)^{1/p}\leq c_p~e^{-\lambda(t-u)}~\left(
 (1+\vert x\vert)~(u-t_n)+(u-t_n)^{1/2}\right)~(u-t_n).
\end{array}
$$
\cqfd

\begin{lem}\label{int-lem1}
For any $t_n\leq u\leq t_{n+1}\leq t$ and $p\geq 1$ we have
$$
\EE\left(\vert\Ja_{t_n,u}(F)(x)\vert^{p}\right)^{1/p}\leq 
c_p \left(1+\vert x\vert\right)~e^{-\lambda(t-u)}~\sqrt{u-t_n}
$$
as well as
$$
\EE\left(\vert \Ka_{t_n,u}(F)(x)\vert^{p}\right)^{1/p}\vee \EE\left(\vert \La_{t_n,u}(F)(x)\vert^{p}\right)^{1/p}\\
\\
\displaystyle\leq c_p ~\left(1+\vert x\vert\right)~e^{-\lambda(t-u)}~(u-t_n). 
$$
\end{lem}
\proof
Using the generalized Minkowski inequality, observe that
\begin{eqnarray*}
\EE\left(\left( \Ja_{t_n,u}(F)(x)\right)^{2p}\right)^{1/p}&\leq& c_p ~\EE\left(
\left(\int_{t_n}^u \left(\frac{u-s}{t_{n+1}-t_n}\right)^2~
\left(\sigma
\partial F_{s,t}\right)(X_{t_n,s}(x))^2~ds\right)^p\right)^{1/p}\\
\\
&\leq &c_p~\int_{t_n}^u \left(\frac{u-s}{t_{n+1}-t_n}\right)^2~
\left(\EE\left(\left(\sigma
\partial F_{s,t}\right)(X_{t_n,s}(x))^{2p}\right)^{1/(2p)}\right)^2~ds\\
\\
&\leq & c_p \left(1+\vert x\vert^2\right)~e^{-2\lambda(t-u)}~\int_{t_n}^u \left(\frac{u-s}{t_{n+1}-t_n}\right)^2~e^{-2\lambda(u-s)}
~ds.
\end{eqnarray*}
Then one has
$$
\begin{array}{l}
\displaystyle\EE\left(\vert \Ka_{t_n,u}(F)(x)\vert^{p}\right)^{1/p}\vee \EE\left(\vert \La_{t_n,u}(F)(x)\vert^{p}\right)^{1/p}\\
\\
\displaystyle\leq c_p ~\left(1+\vert x\vert\right)~e^{-\lambda(t-u)}~ ~\int_{t_n}^{u}~\frac{u-s}{t_{n+1}-t_n}~e^{-\lambda(u-s)}~~ds.
\end{array}
$$
\cqfd

\begin{lem}\label{lem-delta-M}
For any $t_n\leq u\leq t_{n+1}\leq t$  we have
$$
\EE\left(\Ma_{t_n,u}(F)(x)^{2}\right)^{1/2} \leq c ~ \left(1+\vert x\vert\right)~e^{-\lambda(t-u)}~(u-t_n)
$$
and
$$
\begin{array}{l}
\displaystyle \EE\left(\Na_{t_n,u}(F)(x)^{2}\right)^{1/2} \leq c~ \sqrt{\epsilon}~\left(1+\vert x\vert\right)~e^{-\lambda(t-u)} (u-t_n).\end{array}$$
\end{lem}
\proof
Using Lemma~\ref{int-lem1} we check that
$$
\begin{array}{l}
\displaystyle\EE\left(\left( \Ma_{t_n,u}(F)(x)\right)^{2}\right)^{1/2}\\
\\
\displaystyle\leq 
\EE\left(\left(B_{u}-B_{t_n}\right)^{4}\right)^{1/4}\EE\left(\left(\Ja_{t_n,u}(F)(x)\right)^{4}\right)^{1/4}+\EE\left(\left(  \Ka_{t_n,u}(F)(x)\right)^{2}\right)^{1/2}\\
\\
\displaystyle\leq c ~ \left(1+\vert x\vert\right)~e^{-\lambda(t-u)}~(u-t_n).
\end{array}$$
Now for any $t_n\leq u\leq t_{n+1}\leq t$ we have
$$
\begin{array}{l}
 \EE\left(\Na_{t_n,u}(F)(x)^{2}\right)^{1/2}\\
 \\
\displaystyle \leq   \EE\left(\left(B_u-B_{t_n}\right)^{4}\right)^{1/4}\EE\left(\left(\La_{t_n,u}(F)(x)\right)^{4}\right)^{1/4}\\
\\
\displaystyle\hskip3cm+\int_{t_n}^{u}~\frac{u-s}{t_{n+1}-t_n}~ \EE\left(
\left(\left(\sigma
\partial F_{t_n,t}\right)(x)-\left(\sigma
\partial F_{s,t}\right)(X_{t_n,s}(x))\right)^{2}\right)^{1/2}ds.
\end{array}$$
Combining \eqref{contiuity-cond} with Lemma~\ref{int-lem1} we find that
$$
\begin{array}{l}
 \EE\left(\Na_{t_n,u}(F)(x)^{2}\right)^{1/2}\\
 \\
\displaystyle \leq c~  (u-t_n)^{3/2}~~\left(1+\vert x\vert\right)~e^{-\lambda(t-u)} \\
\\
\displaystyle\hskip.3cm+e^{-\lambda(t-u)} \int_{t_n}^{u}~\frac{u-s}{t_{n+1}-t_n}~e^{-\lambda(u-s)}~~\left(
 (1+\vert x\vert)~(s-t_n)+\sqrt{s-t_n}\right)~ds
\end{array}$$
from which we check that
$$
\begin{array}{l}
\displaystyle \EE\left(\Na_{t_n,u}(F)(x)^{2}\right)^{1/2} \leq c~  (u-t_n)^{3/2}~~\left(1+\vert x\vert\right)~e^{-\lambda(t-u)} \\
\\
\displaystyle\hskip3cm+e^{-\lambda(t-u)}\left(  (1+\vert x\vert)(u-t_n)^2+\sqrt{t_{n+1}-t_n}~(u-t_n)\right).
\end{array}$$
This ends the proof of the lemma.
\cqfd

\subsection{Proof of Theorem~\ref{theo-unif-int-strato}} \label{theo-unif-int-strato-proof}
For any $t_n\leq t\leq t_{n+1}$ we set
$$
M_{t}(F)(x):=\sum_{0\leq k<n}\Ma_{t_k,t_{k+1}}(F)\left(X_{t_k}(x)\right)+\Ma_{t_n,t}(F)\left(X_{t_n}(x)\right).
$$
For any $l<k$ we have
$$
\EE\left(
\Ma_{t_k,t_{k+1}}(F)\left(X_{t_l,t_k}(x)\right)\right)=0.
$$
This implies that
$$
\begin{array}{l}
\EE\left(\Ma_{t_l,t_{l+1}}(F)\left(X_{t_l}(x)\right)\Ma_{t_k,t_{k+1}}(F)\left(X_{t_l,t_k}(X_{t_l}(x))\right)\right)=0
\end{array}
$$
and therefore
$$
\EE(M_{t}(F)(x)^2):=\sum_{0\leq k<n}\EE\left(\Ma_{t_k,t_{k+1}}(F)\left(X_{t_k}(x)\right)^2\right)+\EE\left(\Ma_{t_n,t}(F)\left(X_{t_n}(x)\right)^2\right).
$$
Applying Lemma~\ref{lem-delta-M} yields
$$
\begin{array}{l}
\EE(M_{t}(F)(x)^2)\\
\\
\displaystyle \leq c_1~ \epsilon~\left(1+\vert x\vert^2\right)~\left(\sum_{0\leq k<n}e^{-2\lambda(t-t_{k+1})}~(t_{k+1}-t_k)
+(t_{n+1}-t_n)\right)\\
\\
\displaystyle \leq c_2~ \epsilon~\left(1+\vert x\vert^2\right)~
\left(\epsilon+e^{2\lambda\epsilon}~\int_0^{t_n}~e^{-2\lambda(t-s)}~ds
\right)\leq 
c_3~ \epsilon~\left(1+\vert x\vert^2\right)~
\left(\epsilon+e^{2\lambda\epsilon}~e_{2\lambda}(t)
\right).
\end{array}
$$
For any $t_n\leq t\leq t_{n+1}$ we set
$$
I_{t}(F)(x):=\sum_{0\leq k<n}\Ia_{t_k,t_{k+1}}(F)\left(X_{t_k}(x)\right)+\Ia_{t_n,t}(F)\left(X_{t_n}(x)\right).
$$
For any $k<l$ we have
$$
\EE\left(
\Ia_{t_k,t_{k+1}}(F)\left(X_{t_k}(x)\right)~|~X_{t_k}(x)\right)=0.
$$
Applying Lemma~\ref{lem-delta-M} and arguing as above we check that
$$
\begin{array}{l}
\EE(I_{t}(F)(x)^2)\\
\\
\displaystyle \leq c_1~\epsilon~ \left(1+
 \vert x\vert\epsilon^{1/2}\right)^2~\left(\sum_{0\leq k<n}
~e^{-2\lambda(t-t_{k+1})}~~(t_{k+1}-t_k)+~(t_{n+1}-t_n)\right)\\
\\
\leq \displaystyle  c_2~\epsilon~ \left(
 1+\vert x\vert~\epsilon^{1/2}\right)^2\left(\epsilon+e^{2\lambda\epsilon}~e_{2\lambda}(t)
\right).
\end{array}
$$

For any $t_n\leq t\leq t_{n+1}$ we set
$$
N_{t}(F)(x):=\sum_{0\leq k<n}\Na_{t_k,t_{k+1}}(F)\left(X_{t_k}(x)\right)+\Na_{t_n,t}(F)\left(X_{t_n}(x)\right).
$$
Once again using Lemma~\ref{lem-delta-M} we have
$$
\begin{array}{l}
\displaystyle
\EE\left(N_{t}(F)(x)^2\right)^{1/2}\\
\\
\displaystyle\leq c~ \sqrt{\epsilon}~\left(1+\vert x\vert\right)~\left(\sum_{0\leq k<n}~e^{-\lambda(t-t_{k+1})} (t_{k+1}-t_k)+ (t_{n+1}-t_n)\right)\\
\\
\displaystyle\leq c~ \sqrt{\epsilon}~\left(1+\vert x\vert\right)~\left(\epsilon+e^{\lambda\epsilon}~e_{\lambda}(t)
\right).
\end{array}
$$
For any $t_n\leq t\leq t_{n+1}$ we set
$$
D_{t}(F)(x):=\sum_{0\leq k<n}\Da_{t_k,t_{k+1}}(F)\left(X_{t_k}(x)\right)+\Da_{t_n,t}(F)\left(X_{t_n}(x)\right).
$$
Then using Lemma~\ref{lem-delta-M}  and arguing as above for any $p\geq 1$ we have
$$
\EE\left(\vert D_{t}(F)(x)\vert^p\right)^{1/p}\leq 
c_p~\epsilon^{1/2}~\left(
 1+\vert x\vert~\epsilon^{1/2}\right)\left(\epsilon+e^{\lambda\epsilon}~e_{\lambda}(t)
\right).
$$
Recalling that ${\Oa}_{t_k,t_{k+1}}(F)(x)=0$, for any $t_n\leq t\leq t_{n+1}$ we have 
$$
O_{t}(F)(x):=\sum_{0\leq k<n}\Oa_{t_k,t_{k+1}}(F)\left(X_{t_k}(x)\right)+\Oa_{t_n,t}(F)\left(X_{t_n}(x)\right)=\Oa_{t_n,t}(F)\left(X_{t_n}(x)\right).
$$
In addition, applying Lemma~\ref{int-lem1} we have
$$
\begin{array}{l}
\displaystyle
\EE\left(\vert {\Oa}_{t_n,u}(F)(x)-\EE({\Oa}_{t_n,u}(F)(x))\vert^p\right)^{1/p}\\
\\
\displaystyle\leq c_p~(t_{n+1}-u)^{1/2}~ \left(1+\vert x\vert\right)~e^{-\lambda(t-u)}~\sqrt{u-t_n}\\
\\
\displaystyle+(1+\vert x\vert)~
\EE\left(\left\vert \left(\frac{u-t_n}{t_{n+1}-t_n}\left(B_{t_{n+1}}-B_{u}\right)~
-\left(B_u-B_{t_n}\right)~\frac{t_{n+1}-u}{t_{n+1}-t_n}\right)\right\vert^p\right)^{1/p}
\end{array}$$
from which we check that
$$
\begin{array}{l}
\displaystyle
\EE\left(\vert O_{t}(F)(x)\vert^p\right)^{1/p}\leq c_p~\left(1+\vert x\vert\right)~\sqrt{\epsilon}.
\end{array}$$

Applying Proposition~\ref{prop-int-over-B-2}, for any $t_n\leq t\leq t_{n+1}$ we have
$$
\begin{array}{l}
\displaystyle\int_{0}^{t}~F_{s,t}(X_{s}(x))~d\overline{B}_s
=\displaystyle\sum_{0\leq k<n}\int_{t_k}^{t_{k+1}}~F_{s,t}\left(X_{t_k,s}\left(X_{t_k}(x)\right)\right)~d\overline{B}_s+\int_{t_n}^tF_{s,t}(X_{t_n,s}\left(X_{t_n}(x)\right))~d\overline{B}_s\\
\\
\displaystyle=\int_{0}^{t}~F_{s,t}(X_{s}(x))~dB_s+\frac{1}{2}
\int_{0}^{t}~\left(\sigma\partial F_{s,t}\right)(X_{s}(x))~ds+R_{t}(F)(x)
\end{array}
$$
with the remainder term
$$
\begin{array}{l}
\displaystyle R_{t}(F):=M_t(F)+N_t(F)+I_t(F)+D_t(F)+O_t(F)
\end{array}
$$
such that
$$
\EE(R_t(F)(x)^2)^{1/2}\leq c~\sqrt{\epsilon}~(1+e_{\lambda}(t)\vee e_{2\lambda}(t)^{1/2} )~(1+\vert x\vert).
$$
The proof of Theorem~\ref{theo-unif-int-strato} is thus complete.\cqfd

\subsubsection*{Acknowledgements}

AJ \& HR are supported by KAUST baseline funding. In addition, this work was also supported by the Innovation and Talent Base for Digital Technology and Finance (B21038) and `the Fundamental Research Funds for the Central Universities', Zhongnan University of Economics and Law (2722023EJ002).

\section*{Appendix}
\subsection*{Proof of \eqref{LG-stable} and \eqref{LG-unstable}}
We have
$$
X_{t_n}(x)-\overline{X}_{t_n}(x)=\sum_{1\leq k\leq n}\int_{t_{k-1}}^{t_{k}}\left(e^{a(t_n-s)}-\frac{1}{\epsilon}\int_{t_{k-1}}^{t_{k}}~e^{a(t_n-u)}du\right)~dB_s.
$$
This already shows that the smooth approximating model is unbiased and its variance is given by
$$
\begin{array}{l}
\EE\left((X_{t_n}(x)-\overline{X}_{t_n}(x))^2\right)\\
\\
\displaystyle=\sum_{1\leq k\leq n}~e^{2a(t_n-t_k)}~\left(\frac{1}{\epsilon}~\int_{t_{k-1}}^{t_k}\left(e^{a(t_k-s)}-\frac{1}{\epsilon}~\int_{t_{k-1}}^{t_k}e^{a(t_k-u)}~ du\right)^2 ds\right)~
(t_k-t_{k-1}).
\end{array}
$$
We note that
$$
\begin{array}{rcl}
\displaystyle\frac{1}{\epsilon}~\int_{0}^{\epsilon}\left(e^{a s}-\frac{1}{\epsilon}~\int_{0}^{\epsilon}e^{au}~ du\right)^2~ds&=&\displaystyle\frac{1}{2}~\frac{1}{\epsilon^2}\int_{[0,\epsilon]^2}~(e^{as}-e^{au})^2~duds\\
\\
&=&\displaystyle a^2~\frac{1}{2\epsilon^2}\int_{[0,\epsilon]^2}(s-u)^2~\left(\int_0^1e^{a(\delta s+(1-\delta) u)} d\delta\right)^2~duds.
\end{array}
$$
In addition, we have
$$
\frac{1}{2\epsilon^2}\int_{[0,\epsilon]^2}~(s-u)^2 ds du=\frac{1}{\epsilon}\int_{0}^{\epsilon}~s^2 ds-\left(\frac{1}{\epsilon}\int_{0}^{\epsilon}~u du\right)^2=\epsilon^2/12.
$$
Now we consider the following two cases:
\begin{itemize}
\item When $a<0$ we have
$$
\begin{array}{l}
\EE\left((X_{t_n}(x)-\overline{X}_{t_n}(x))^2\right)\\
\\
\displaystyle\leq \frac{(a\epsilon)^2}{12}~\sum_{1\leq k\leq n}~\int_{t_k}^{t_{k+1}}~e^{2a(t_n-s)}~ds\\
\\
\displaystyle= a~\frac{\epsilon^2}{24}~\int_{t_1}^{t_{n+1}}~(2a)~e^{2a(t_n-s)}~ds=\vert a\vert~e^{2\vert a\vert \epsilon}~\frac{\epsilon^2}{24}~\left(1-e^{2at_n}\right).
\end{array}
$$
This yields the uniform estimate \eqref{LG-stable}.

\item When $a>0$ we have
$$
\frac{\epsilon^2}{12}~a^2\leq \frac{1}{\epsilon}~\int_{0}^{\epsilon}\left(e^{a s}-\frac{1}{\epsilon}~\int_{0}^{\epsilon}e^{au}~ du\right)^2~ds\leq \frac{\epsilon^2}{12}~a^2~e^{2a\epsilon}
$$

$$
\begin{array}{l}
\EE\left((X_{t_n}(x)-\overline{X}_{t_n}(x))^2\right)\\
\\
\displaystyle\leq \frac{\epsilon^2}{12}~a^2~e^{2a\epsilon}
~\sum_{1\leq k\leq n}~\int_{t_{k-1}}^{t_k}~e^{2a(t_n-s)}~ds=
\epsilon^2~\frac{a}{24}~e^{2a\epsilon}~\left(e^{2at_n}-1\right)
\end{array}
$$
$$
\begin{array}{l}
\EE\left((X_{t_n}(x)-\overline{X}_{t_n}(x))^2\right)\\
\\
\displaystyle\geq \frac{\epsilon^2}{12}~a^2~
~\sum_{1\leq k\leq n}~\int_{t_{k}}^{t_{k+1}}~e^{2a(t_n-s)}~ds=\frac{\epsilon^2}{24}~a~\int_{t_1}^{t_{n+1}}~2a~e^{2a(t_n-s)}~ds=\epsilon^2~\frac{a}{24}~e^{-2a\epsilon}~\left(e^{2at_n}-1\right).
\end{array}
$$
This yields the estimate \eqref{LG-unstable}.
\end{itemize}
The proof of  \eqref{LG-stable} and \eqref{LG-unstable} is now completed.\cqfd

\subsection*{Proof of \eqref{lamperti}}\label{lamperti-proof}
$$
\theta(x):=\int_{0}^x~\frac{1}{\sigma(y)}~dy\Longrightarrow \partial \theta=
\frac{1}{\sigma}\quad \mbox{\rm and}\quad
\partial^2 \theta=-
\frac{1}{\sigma^2}~\partial \sigma.
$$
Now for any $x_1\leq x_2$ we have
$$
(x_2-x_1)/\sigma_+\leq\theta(x_2)-\theta(x_1)=\int_{x_1}^{x_2}~\frac{1}{\sigma(y)}~dy\leq  (x_2-x_1)/\sigma_-.
$$
On the other hand, applying the It\^o formula we have
$$
d\theta(X_t(x))=\frac{1}{\sigma(X_t(x))}~\left(b_{\sigma}\left(X_{t}(x)\right)~dt+~\sigma\left(X_{t}(x)\right)
    dB_t\right)-\frac{1}{2\sigma(X_t(x))^2}~\partial \sigma \left(X_{t}(x)\right)~\sigma\left(X_{t}(x)\right)^2~dt.
$$
This implies that
$$
d\theta(X_t(x))=b_{/\sigma}\left(X_{t}(x)\right)~dt+
    dB_t\quad \mbox{\rm with}\quad b_{/\sigma}:=b/\sigma.
$$
In the same vein, we have
$$
d\theta(\overline{X}_t(x))=b_{/\sigma}\left(\overline{X}_{t}(x)\right)~dt+
    d\overline{B}_t.
    $$
 Rewritten in terms of the processes
 $$
 X_t^{\theta}(x):=\theta(X_t(x))\quad\mbox{\rm and}\quad
  \overline{X}_t^{\theta}(x):=\theta(\overline{X}_t(x))
 $$       
we have proved that
$$
dX^{\theta}_t(x)=b^{\theta}\left(X^{\theta}_{t}(x)\right)~dt+dB_t
\quad\mbox{\rm and}\quad  d\overline{X}_t^{\theta}(x)=b^{\theta}\left(\overline{X}^{\theta}_{t}(x)\right)~dt+d\overline{B}_t
$$
with the drift function
$$
 b^{\theta}:=b_{/\sigma}\circ\theta^{-1}=\left(\left(\partial\theta\right)\circ\theta^{-1}\right)~\left(b\circ\theta^{-1}\right).
$$
\cqfd

\subsection*{Proof of Lemma~\ref{lem-Z-exp}}\label{lem-Z-exp-proof}

The gradient $\nabla Z_{s,t}(x)$ of the  flow $Z_{s,t}(x)$ is given by  the $(r\times r)$-matrix
$$
d  \,\nabla Z_{s,t}(x)=\nabla Z_{s,t}(x)~\nabla b\left(Z_{s,t}(x)\right)~dt\quad \mbox{\rm with}\quad \nabla Z_{s,s}(x)=I.
$$
After some calculations, one has
$$
 \begin{array}{rcl}
\displaystyle d  \,\left[\nabla Z_{s,t}(x) \,\nabla Z_{s,t}(x)^{\prime}\right]
&=&\nabla Z_{s,t}(x) ~\left(\nabla b\left(Z_{s,t}(x)\right)+\nabla b\left(Z_{s,t}(x)\right)^{\prime}\right)~\nabla Z_{s,t}(x)^{\prime}~dt\\
&&\\
&\leq &-2\lambda_b~\left[\nabla Z_{s,t}(x) \,\nabla Z_{s,t}(x)^{\prime}\right]~dt.
\end{array} 
$$
Taking the trace, we check that
$$
\Vert \nabla Z_{s,t}\Vert_{\tiny Frob}:=
\sup_{x\in\RR^r}\Vert \nabla Z_{s,t}(x)\Vert_{\tiny Frob}\leq e^{-\lambda_b (t-s)}.
$$
The Taylor expansion
$$
\begin{array}{l}
\displaystyle Z_{t}(x)-Z_{t}(y)=\int_0^1~\nabla Z_{t}(\epsilon x+(1-\epsilon)y)^{\prime}(x-y)~d\epsilon\\
\\
\displaystyle\Longrightarrow \Vert Z_{t}(x)-Z_{t}(y)\Vert\leq  c_1~e^{-\lambda_b t}~\Vert x-y\Vert
\end{array}
$$
ensures that $Z_t(x)$ converges exponentially fast  towards  a fixed point $x_0=Z_t(x_0)\in\RR^r$  as $t\rightarrow\infty$. This yields the uniform estimate
$$
\Vert Z_{t}(x)\Vert\leq \Vert x_0\Vert+c_1~e^{-\lambda_b t}~\Vert x-x_0\Vert\leq c_2~(1\vee \Vert x_0\Vert)
( 1+\Vert x\Vert).
$$

Now we have
$$
 \begin{array}{l}
d \, \nabla^2 Z_{s,t}(x)\\
\\
=\left[\left[\nabla Z_{s,t}(x)\otimes \nabla Z_{s,t}(x)\right]\nabla^2b(Z_{s,t}(x))+\nabla^2 Z_{s,t}(x)\nabla b(Z_{s,t}(x))\right]dt\quad \mbox{\rm with}\quad \nabla^2 Z_{s,s}(x)=0.
\end{array}$$
Then one can show that
$$
 \begin{array}{l}
\displaystyle \partial_t\, \left[\nabla^2 Z_{s,t}(x)\nabla^2 Z_{s,t}(x)^{\prime}\right]\\
\\
\displaystyle=\left[\nabla^2 Z_{s,t}(x)~(\nabla b_t (Z_{s,t}(x)))_{\tiny sym}~\nabla^2 Z_{s,t}(x)^{\prime}\right]+2\left[\left[\nabla Z_{s,t}(x)\otimes \nabla Z_{s,t}(x)\right]~\nabla^2b(Z_{s,t}(x))~\nabla^2 Z_{s,t}(x)^{\prime}\right]_{\tiny sym}.
\end{array}$$
Whenever $$
\Vert\nabla^2b\Vert_{\tiny Frob}:=
\sup_{x\in\RR^r}\Vert\nabla^2b(x)\Vert_{\tiny Frob}<\infty$$  taking the trace in the above display we check that
$$
 \partial_t\, \Vert\nabla^2 Z_{s,t}(x)\Vert_{\tiny Frob}^2\leq -2\lambda_b~\Vert\nabla^2 Z_{s,t}(x)\Vert_{\tiny Frob}^2+2\Vert\nabla^2b\Vert_{\tiny Frob}~\Vert\nabla Z_{s,t}(x)\Vert_{\tiny Frob}^2~\Vert\nabla^2 Z_{s,t}(x)\Vert_{\tiny Frob}.
$$
This yields the estimate
$$
 \partial_t\, \Vert\nabla^2 Z_{s,t}(x)\Vert_{\tiny Frob}\leq -\lambda_b~\Vert\nabla^2 Z_{s,t}(x)\Vert_{\tiny Frob}+\Vert\nabla^2b\Vert_{\tiny Frob}~\Vert\nabla Z_{s,t}(x)\Vert_{\tiny Frob}^2.
$$
Now we have
$$
\Vert\nabla^2 Z_{s,t}(x)\Vert_{\tiny Frob}\leq \Vert\nabla^2b\Vert_{\tiny Frob}~e^{-\lambda_b(t-s)}~\int_s^t~e^{\lambda_b(u-s)}~\Vert\nabla Z_{s,u}(x)\Vert_{\tiny Frob}^2~du
$$
from which we conclude that
$$
\Vert\nabla^2 Z_{s,t}\Vert_{\tiny Frob}:=\sup_{x\in\RR^r}
\Vert\nabla^2 Z_{s,t}(x)\Vert_{\tiny Frob}\leq \frac{r}{\lambda_b}~\Vert\nabla^2b\Vert_{\tiny Frob}~e^{-\lambda_b(t-s)}.
$$
\cqfd

\subsection*{Proof of Proposition~\ref{interpol}}\label{interpol-proof}

For any $u\leq t$ and any function $f$ on $\RR^r$ we have
the backward formula
$$
d_u (f\circ Z_{u,t})(x)=-\nabla (f\circ Z_{u,t})(x)^{\prime}~b(x)~du
$$
from which we have that
 \begin{eqnarray*}
d_u\left( (f\circ Z_{u,t})\circ \overline{X}_{s,u}\right)
&=&\sum_{1\leq j\leq \overline{r}}~\partial_{\sigma_j}(f\circ Z_{u,t})(\overline{X}_{s,u})~d\overline{B}^j_u.\end{eqnarray*}
This yields the interpolation formula
 \begin{eqnarray*}
f(\overline{X}_{s,t})&=&f(Z_{s,t})+\sum_{1\leq j\leq \overline{r}}~\int_s^t~\partial_{\sigma_j}(f\circ Z_{u,t})(\overline{X}_{s,u})~d\overline{B}^j_u.
 \end{eqnarray*}

Applying the above formula to $\partial_{\sigma_j}(f\circ Z_{u,t})$ we also have
 \begin{eqnarray*}
\partial_{\sigma_j}(f\circ Z_{u,t})(\overline{X}_{s,u})&=&\partial_{\sigma_j}(f\circ Z_{u,t})(Z_{s,u})+\sum_{1\leq i\leq \overline{r}}\int_s^u~\partial_{\sigma_i}\left(\partial_{\sigma_j}(f\circ Z_{u,t})\circ Z_{u,v}\right)(\overline{X}_{s,v})~d\overline{B}^i_v
 \end{eqnarray*}
from which we conclude.
\cqfd

\subsection*{Proof of \eqref{hyp-mom}}\label{hyp-mom-proof}
For any time horizon $s\geq 0$ denote by
  $Z^{\sigma}_{s,t}(x)$  the (deterministic) flow defined for any $t\in [s,\infty[$ and any starting point $Z^{\sigma}_{s,s}(x)=x\in\RR^r$ by the ordinary differential equation
$$
dZ^{\sigma}_{s,t}(x)=b_{\sigma}\left(Z^{\sigma}_{s,t}(x)\right)~dt.
$$
Whenever  $(H_{b_{\sigma}})$ is satisfied, Lemma~\ref{lem-Z-exp} ensures that
$$
\Vert\nabla Z^{\sigma}_{s,t}(x)\Vert\vee \Vert\nabla^2 Z^{\sigma}_{s,t}(x)\Vert\leq c~e^{-\lambda_{b_{\sigma}}(t-s)}
\quad\mbox{and}\quad 
\Vert Z^{\sigma}_{s,t}(x)\Vert\leq c~(1+\Vert x\Vert).
$$
This implies $Z^{\sigma}_t(x)$ converges exponentially fast  towards  a fixed point $x_0=Z^{\sigma}_t(x_0)\in\RR^r$  as $t\rightarrow\infty$. Therefore we have the uniform estimate
$$
\Vert Z^{\sigma}_{t}(x)\Vert\leq \Vert x_0\Vert+c_1~e^{-\lambda_{b_{\sigma}} t}~\Vert x-x_0\Vert\leq c_2~(1\vee \Vert x_0\Vert)
(1+\Vert x\Vert).
$$
Applying \cite[Theorem 1.2]{dmsumeet} we have
the interpolation formula
$$
X_{s,t}(x)-Z^{\sigma}_{s,t}(x)=\int_s^t\nabla_{\sigma} Z^{\sigma}_{u,t}(X_{s,u}(x))^{\prime} ~dB_u+\frac{1}{2}~\int_s^t
 \nabla^2 Z^{\sigma}_{u,t}(X_{s,u}(x))^{\prime}~a(X_{s,u}(x))~du
$$
from which we check that
$$
\sup_{t\geq 0}\EE\left(\Vert X_{t}(x)-Z^{\sigma}_{t}(x)\Vert^p \right)^{1/p}\leq c_p~(1+\Vert x\Vert)
$$
The proof of \eqref{hyp-mom} is now easily completed. \cqfd

\subsection*{Proof of \eqref{fin-s}}\label{fin-s-proof}
For a given uniformly bounded Lipschitz function $f$ on $\RR$, we set
$$
f_{\sigma}:=f/\sigma\quad \mbox{\rm and}\quad g_{\sigma}:=b f_{\sigma}.
$$
We have that
$$
F(x):=\int_0^x f_{\sigma}(y)~ dy\Longrightarrow \partial F=f_{\sigma}.
$$
Applying It\^o formula, we have
\begin{eqnarray*}
dF(X_t(x))&=&f_{\sigma}(X_t(x))~dX_t(x)+\frac{1}{2}~\sigma^2(X_t(x))~\partial f_{\sigma}(X_t(x))~dt\\
dF(\overline{X}_t(x))&=&f_{\sigma}(\overline{X}_t(x))~d\overline{X}_t(x).
\end{eqnarray*}
This implies that
$$
\begin{array}{l}
\displaystyle F(X_t(x))-F(x)-\int_0^t~f_{\sigma}(X_s(x))~b_{\sigma}(X_s(x))~ds\\
\\
\displaystyle=\int_0^t f_{\sigma}(X_s(x))~\sigma(X_s(x))~dB_s+\frac{1}{2}~\int_0^t \sigma^2(X_s(x))~\partial f_{\sigma}(X_s(x))~ds
\end{array}$$
and
$$
F(\overline{X}_t(x))-F(x)-\int_0^t~f_{\sigma}(\overline{X}_s(x))~b(\overline{X}_s(x))~ds=\int_0^t f_{\sigma}(\overline{X}_s(x))~\sigma(\overline{X}_s(x))~d\overline{B}_s.
$$
Letting $g_{\sigma}:=b f_{\sigma}$ we conclude that
$$
\begin{array}{l}
\displaystyle F(\overline{X}_t(x))-F(X_t(x))+\int_0^t\left(g_{\sigma}(\overline{X}_s(x))-g_{\sigma}(X_s(x))\right)~ds\\
\\
\displaystyle=\int_0^t f(\overline{X}_s(x))~d\overline{B}_s
-\left(\int_0^t f(X_s(x))~dB_s+\frac{1}{2}\int_0^t \sigma(X_s(x))~\partial f(X_s(x))~ds
\right).
\end{array}
$$
Whenever $f_{\sigma}$ is uniformly bounded and $g_{\sigma}:=b f_{\sigma}$ is Lipschitz, we have
$$
\EE\left(\left(F(\overline{X}_t(x))-F(X_t(x))
\right)^{2}\right)\leq c_1~\epsilon
$$
and
$$
\EE\left(\left(\int_0^t\left(g(X_s(x))-g(\overline{X}_s(x))\right)~ds
\right)^{2}\right)^{1/2}\leq \int_0^t
\EE\left(\left(g(X_s(x))-g(\overline{X}_s(x))
\right)^{2}\right)^{1/2}~ds\leq~t~c_2~\sqrt{\epsilon}.
$$
Then one can obtain that
$$
\EE\left(\left(\int_0^t f(\overline{X}_s(x))~d\overline{B}_s
-\left(\int_0^t f(X_s(x))~dB_s+\frac{1}{2}\int_0^t \sigma(X_s(x))~\partial f(X_s(x))~ds
\right)
\right)^{2}\right)^{1/2}\leq c~(1+t)~\sqrt{\epsilon}
$$
which completes the proof of~\eqref{fin-s}.\cqfd

\end{document}